%% file: spacedepME.tex
\documentclass[fleqn,12pt]{article}

\usepackage[latin1]{inputenc}
\usepackage{graphicx}
\usepackage{subfigure}
\usepackage{amsmath}
\usepackage{amssymb}
\usepackage{here} 
\usepackage{a4}
\usepackage{natbib}

\def\faty{\mathbf{y}}

\def\fate{\mathbf{e}}
\def\fatf{\mathbf{f}}
\def\fatg{\mathbf{g}}
\def\fatm{\mathbf{m}}
\def\fatn{\mathbf{n}}

\def\fatx{\mathbf{x}}

\def\fatphi{\boldsymbol{\phi}}

\def\fatomega{\boldsymbol{\omega}}

\def\tri{\mathcal{C}}

\def\calD{\mathcal{D}}
\def\calM{\mathcal{M}}

\newcommand{\ba}[1] {\left( \begin{array}{#1}}
\newcommand{\ea}{\end{array} \right) }
\newcommand{\ordo}[1]{{\cal O} \left( #1 \right)}
\newcommand{\natop}[2]{\genfrac{}{}{0pt}{1}{#1}{#2}}

\title{{\Large 
  SIMULATION OF STOCHASTIC
  REACTION-DIFFUSION PROCESSES ON
  UNSTRUCTURED MESHES }}

\author{{\small 
  STEFAN ENGBLOM$^{\mbox{\tiny 1}}$, LARS FERM$^{\mbox{\tiny 1}}$} \\ 
  {\small
  ANDREAS HELLANDER$^{\mbox{\tiny 1}}$, PER LÖTSTEDT$^{\mbox{\tiny 1}}$}
\thanks{Financial support has been obtained from the Swedish 
Foundation for Strategic Research and the Swedish National Graduate
School in Mathematics and Computing. Corresponding author: Per
Lötstedt, address as above, telephone +46-18-4712972, fax
+46-18-523049.}}

\date{\today}

\numberwithin{equation}{section}
\numberwithin{table}{section}
\numberwithin{figure}{section}

\newtheorem{theorem}{Theorem}[section]
\newtheorem{assumption}[theorem]{Assumption}

\newtheorem{proposition}[theorem]{Proposition}

\newenvironment{proof}
{\noindent \textit{Proof.}}
{\hspace{1pt} $\square$ \vspace{3mm}}

\newenvironment{remark}
{\noindent \textit{Remark.}}
{\hspace{1pt} $\square$ \vspace{3mm}}

\begin{document}

\maketitle

\vspace{-10pt}

\begin{center}
{\footnotesize\em 
$^{\mbox{\tiny\rm 1}}$Division of Scientific Computing,
Department of Information Technology \\
Uppsala University, P.~O.~Box 337, SE-75105 Uppsala, Sweden \\
emails: {\tt stefane, ferm, andreas.hellander, perl@it.uu.se} \\
[3pt]}
\end{center}


\begin{abstract}

Stochastic chemical systems with diffusion are modeled with a
reaction-diffusion master equation. On a macroscopic level, the
governing equation is a reaction-diffusion equation for the averages
of the chemical species. On a mesoscopic level, the master equation
for a well stirred chemical system is combined with Brownian motion in
space to obtain the reaction-diffusion master equation. The space is
covered by an unstructured mesh and the diffusion coefficients on the
mesoscale are obtained from a finite element discretization of the
Laplace operator on the macroscale. The resulting method is a flexible
hybrid algorithm in that the diffusion can be handled either on the
meso- or on the macroscale level. The accuracy and the efficiency of
the method are illustrated in three numerical examples inspired by
molecular biology.

\vspace{0.5cm}

\noindent
{\bf Keywords}: master equation, diffusion, reaction rate equations,
hybrid method, finite element method, unstructured mesh.

\vspace{0.5cm}

\noindent
{\bf AMS subject classification}: 65C40, 65C05, 65M60, 60H35.
%
%
%

\end{abstract}


\input nomen
\input intro
\input eq
\input diff
\input mom
\input hybrid

\input res
\input concl

\bibliographystyle{plain}
\bibliography{efhl}

\end{document}

%% file: nomen.tex
\section*{Abbreviations}
\label{sec:nomen}

\begin{table}[H]
\begin{center}
\begin{tabular}{ll}
\hline
Abbreviation & Expanded form \\
\hline
mRNA & messenger ribonucleic acid\\
CME & chemical master equation\\
PDF & probability density function \\
SSA & stochastic simulation algorithm\\
RRE & reaction rate equations\\
ODE & ordinary differential equation\\
RDME & reaction-diffusion master equation\\
NRM & next reaction method\\
NSM & next subvolume method\\
BD & Brownian dynamics\\
PDE & partial differential equation\\
RDE & reaction-diffusion equation\\
$n$D & $n$ space dimensions\\
FEM & finite element method\\
SPDE & stochastic PDE\\
FVM & finite volume method\\
\hline
\end{tabular}
\end{center}
\caption{Abbreviations in the paper in order of appearance.}
\label{tab:abbr}
\end{table}

%% file: intro.tex
\section{Introduction}
\label{sec:intro}

Intrinsic noise in biochemical networks can have a large impact on the
macroscopic behavior of biological cells \cite{ELSS, LoH, MAA2, PBE,
ROS}. An example is the regulation of the transcription of genes to
messenger RNA (mRNA) where genes are present in one or two copies and
the copy number of mRNA is small. The facts that the copy number is a
small nonnegative integer and that there is a probability that a
certain reaction will occur when two molecules meet make a discrete,
stochastic description of the system necessary.

The state of the system is the number of molecules of each
participating species. Usually, the assumption is that the system is
well stirred such that there is no spatial dependence of the
distribution of the chemical species. Then the chemical master
equation (CME) is the governing equation for the probability density
function (PDF) of the state of the chemical system \cite{Gardiner,
VanKampen}. A trajectory of the biochemical system is simulated by
randomly choosing a reaction and then updating the state vector in
Gillespie's Stochastic Simulation Algorithm (SSA) \cite{gillespie},
further developed in \cite{multiscaleSSA, nestedSSA, GibsonBruck,
haseltine_HSSA} to improve the efficiency and to handle systems with
slow and fast time scales. The concentrations of the species at a
macroscopic level are often approximated by the reaction rate
equations (RRE) which is a deterministic system of nonlinear ordinary
differential equations (ODEs). This approach works well only when the
number of molecules is large, a condition which is often violated
inside living cells \cite{guptasarama}.

There are many biochemical systems where the spatial inhomogeneity of
the species cannot be neglected. Such systems are no longer well
stirred since the transport of the molecules through the solvent is
slow compared to typical reaction times \cite{KrE} or since some
reactions are strongly localized. The correlation length, i.e.~the
length scale on which the system can be regarded as spatially
homogeneous, is now much shorter. Examples are found in \cite{DRKB,
DuHo, ElE, FaE, Met} where both the stochastic properties and the
spatial distribution are necessary to explain experimental data. If
the diffusion at a molecular level is treated as a special set of
reactions in the CME, then we arrive at the reaction-diffusion master
equation (RDME) \cite[Ch.~8]{Gardiner}, \cite[Ch.~XIV]{VanKampen},
\cite{MMH}. This is an equation for the time evolution of the PDF of
the state of the system in the same manner as the CME but the
dimensionality of the problem is much higher.

In the RDME, the geometry is partitioned into computational cells and
in each cell there is a stochastic variable representing the number of
molecules of each species. If the number of species is $N$ and the
number of cells is $K$ then the RDME is an equation for the scalar PDF
in $KN$ dimensions and time. For example, for a small problem with
five species and $100$ nodes in a mesh, the dimensionality is $500$
and a direct solution of the RDME is obviously out of the
question. The only feasible way is to generate samples from the RDME
and collect statistics in a Monte Carlo fashion. Examples where the
SSA has been applied to reaction-diffusion systems in one space
dimension can be found in \cite{Ber, StL}. An efficient version of the
SSA is the \textit{next reaction method} (NRM) \cite{GibsonBruck}. An
implementation of the NRM, specially developed for diffusive systems,
is the \textit{next subvolume method} (NSM) \cite{ElE}, implemented in
the freely available computer software MesoRD
\cite{HFE}.

More accurate modeling may be necessary if the number of molecules is
very low. In Brownian dynamics (BD), the separate paths of single
molecules are tracked and they may react with other molecules in their
vicinity \cite{AnB, RKDB, ZoR}. The reaction and diffusion of the
particles are simulated using the Green's function of the Smoluchowski
and diffusion equation in \cite{ZoR}. The BD approach is possible only
if the total number of molecules is small and becomes inefficient when
there are many nonreactive collisions for each reactive one. The RDME
and BD are compared in \cite{DRKB}.

The corresponding macroscopic equation for the concentrations of the
species is a partial differential equation (PDE): the
reaction-diffusion equation (RDE). This is the equation solved in
applications such as combustion \cite{PoV} and the model is
appropriate when the number of molecules is large and stochastic
effects can be neglected.

In this paper, we develop a method for simulating the RDME on an
unstructured mesh consisting of triangles in two space dimensions (2D)
or tetrahedra in 3D. Unstructured meshes have the advantage of
approximating curved inner and outer boundaries much more accurately
than Cartesian meshes do. The diffusion coefficients at the meso level
are chosen to be consistent with the discretization with the finite
element method (FEM) on the mesh converging to the diffusion operator
at the macro level when the cell size vanishes. With a proper
computational mesh \cite{Geo}, the FEM approximation yields positive
probabilities for a particle to jump into the adjacent cell. The time
integration of the RDME is split according to Strang \cite{Str} into
two parts in a hybrid method: the diffusion and the chemical
reactions. First, a macroscopic diffusive step is taken for a subset
of the species with the diffusion operator at the macro level using
the unstructured primal mesh. Then the stochastic reactions and the
stochastic diffusion for the remaining species are advanced in the
dual cells of the mesh with SSA at the meso level. If all species in
all cells are treated at the meso level, then we exactly sample the
RDME.



Hybrid methods for a well stirred system with a mesoscopic-macroscopic
approximation are found in \cite{haseltine_HSSA, hybridHell}. In
\cite{MoC}, an efficient Monte Carlo method for the diffusion equation
is described. In \cite{IsP}, the mesoscopic diffusion coefficients are
derived from the discretization of the Laplacian on a 2D Cartesian
mesh with an interior boundary. With unstructured meshes, we propose
in this paper to obtain those coefficients from a proper FEM
discretization.


The outline of the paper is as follows. In Section~\ref{sec:eq}, the
RDME is stated and the relation to the RDE is discussed. The diffusion
coefficients are derived from the discrete Laplacian in
Section~\ref{sec:diff}. The first and second moments of the
distribution of molecules in a system without chemical reactions are
derived in Section~\ref{sec:mom}. The hybrid method coupling the meso-
and macroscales and the operator splitting in time are found in
Section~\ref{sec:hybrid}. Three examples in 2D are found in the
section with numerical results. It is shown that the suggested
mesoscopic diffusion yields trajectories converging to the macroscopic
diffusion equation. It is also shown that the hybrid method we propose
accurately samples the RDME at a fraction of the time needed for a
full simulation. Conclusions are drawn in the final section.


%% file: eq.tex
\section{Reaction-Diffusion Master Equation}
\label{sec:eq}

Assume that the computational domain $\Omega$ in space is partitioned
into computational cells $\tri_j,\; j=1,\ldots, K,$ such that the
cells do not overlap and they cover the whole domain
\[
  \tri_i\cap\tri_j=\emptyset, i\neq j,
  \quad {\rm and} \quad 
  \cup_{j=1}^K \tri_j=\Omega.
\]
Furthermore, assume that there are $N$ chemically active species
$X_{ij},\; i=1,\ldots,N$, in the $K$ cells, $j=1,\ldots,K$. The state
of the system is the array $\fatx$ with $N \times K$ components
$x_{ij}$. The $j$th column of $\fatx$ is denoted by $\fatx_{\cdot j}$
and the $i$th row by $\fatx_{i \cdot}$. The non-negative integer
$x_{ij}$ is thus the copy number of species $i$ in cell $j$. The time
dependent state is changed by chemical reactions occurring between the
molecules in the same cell and by diffusion where molecules move to
adjacent cells. In the reactions, the species interact vertically in
the array $\fatx$ and in the diffusion, the interaction is
horizontal. The RDME governs the time evolution of the PDF $p$, where
$p(\fatx,t)$ is the probability to be in state $\fatx$ at time $t$.

\subsection{Chemical reactions}

A reaction $r$ in a cell $j$ is a transition from one state
$\tilde{\fatx}_{\cdot j}$ before the reaction to the state
$\fatx_{\cdot j}=\tilde{\fatx}_{\cdot j}-\fatn_r$ after the
reaction. The state-change vector $\fatn_r$ of a reaction is a vector
with small integer numbers of length $N$ independent of $j$. There is
a reaction probability or \emph{propensity} $w_r$ that reaction $r$
will take place in a cell depending on the state $\tilde{\fatx}_{\cdot
j}$. A chemical reaction in cell $j$ can be written
\begin{equation}
  \tilde{\fatx}_{\cdot j}\xrightarrow{w_r(\tilde{\fatx}_{\cdot j})} 
  \fatx_{\cdot j},\; 
  \tilde{\fatx}_{\cdot j} = \fatx_{\cdot j}+\fatn_r.
  \label{eq:masterreac}
\end{equation}

In a system without diffusion, the PDF for the molecular distribution
in $\tri_j$ satisfies the CME. Let $\fatn_r$ be split into two parts
\[
  \fatn_r=\fatn_r^+ +\fatn_r^-, \; n_{ri}^+=\max(n_{ri}, 0),
  n_{ri}^-=\min(n_{ri}, 0).
\]
Then the CME for $p$ is, see \cite[Ch.~7]{Gardiner},
\cite[Ch.~V]{VanKampen},
\begin{align}
  \label{eq:defCME}
  \frac{\partial p(\fatx, t)}{\partial t} &=
  \calM p (\fatx, t) \equiv                                               \\
  \nonumber
  &\sum_{j=1}^{K}
  \sum_{\natop{r = 1}{\fatx_{\cdot j}+\fatn_r^- \ge 0}}^{R}
  w_r(\fatx_{\cdot j}+\fatn_r)
  p(\fatx_{\cdot 1},\ldots,\fatx_{\cdot j}+\fatn_r,
  \ldots,\fatx_{\cdot K}, t)                                              \\
  \nonumber
  -&\sum_{j=1}^{K}
  \sum_{\natop{r = 1}{\fatx_{\cdot j}-\fatn_r^+ \ge 0}}^{R}
  w_r(\fatx_{\cdot j})p(\fatx, t),
\end{align}
where the constraints on $\fatx$ are defined elementwise. These
constraints are introduced in order to avoid unfeasible reactions and
will be dropped in the following as is customary.

A simple reaction in $\tri_k$ is
\begin{equation}
  \label{eq:exreac}
   X_{ik}+X_{jk}\rightarrow X_{lk},\; w_1(\fatx_{\cdot k}) = 
  c_{1k} x_{ik}x_{jk},
\end{equation}
where conventionally, we use uppercase letters to denote molecule
\emph{names}, while lowercases are used for counting the number of
molecules of a certain species. The transition vector $\fatn_1$ is
zero except for three components: $n_{1i} = n_{1j} = 1, n_{1l} =
-1$. The propensity $w_1$ has the same form for all cells $k$, whereas
the reaction coefficient $c_{1k}$ scales with the area or volume
$|\tri_k|$ of the cell such that $c_{1k} = \hat{c}_1/|\tri_k|$, where
$\hat{c}_1$ is a constant.

\subsection{Diffusion}


Suppose now that there are no chemical reactions but only diffusion in
the system. Then the mesoscale model of the diffusion of species $i$
from one cell $\tri_{k}$ to another cell $\tri_{j}$ can be written as
a chemical reaction (cf.~(\ref{eq:exreac}))
\begin{equation}
  \label{eq:exdiff}
  X_{ik}\rightarrow X_{ij},\; v_{kj}(\fatx_{i \cdot}) =
  q_{kj} x_{ik}.
\end{equation}
It is understood that $q_{kj}$ is non-zero only for those cells that
are connected and $q_{jj} = 0$. The form of the propensity $v_{kj}$
and the diffusion coefficient $q_{kj}$ are the same for all species
here but $q_{kj}$ may depend on $i$ and can be different for small and
large molecules. We can write
\begin{align}
  \label{eq:diffusion_comps}
  q_{kj} &= \gamma \frac{\hat{q}_{kj}}{h_{kj}^{2}},
\end{align}
where $\gamma$ is the macroscopic diffusion constant, $h_{kj}$ is a
measure of the local length-scale and $\hat{q}_{kj}$ is dimensionless
but still depends on the precise shapes of the cells $\tri_{k}$ and
$\tri{j}$. The interpretation of $q_{kj}$ as the inverse of the
expected value of the \emph{first exit time} for a single molecule
from cell $\tri_{k}$ to $\tri_{j}$ makes it clear why no simple
formula exists except for very regular cells. The molecular movement
is often modeled by the It\^{o}-diffusion
\begin{align}
  \label{eq:Ito}
  d\xi &= \sigma dW_{t},
\end{align}
where $\xi(t)$ is the position and $W_{t}$ is a Wiener-process. The
relation between $\gamma$ in~(\ref{eq:diffusion_comps}) and $\sigma$
in~(\ref{eq:Ito}) is then simply $\gamma = \sigma^{2}/2$.

With this notation, the master equation for the diffusion in
(\ref{eq:exdiff}) can be written in the same manner as the CME in
(\ref{eq:defCME}), see \cite[Ch.~8]{Gardiner},
\cite[Ch.~XIV]{VanKampen}, \cite{MMH},
\begin{align}
  \nonumber
  \frac{\partial p(\fatx, t)}{\partial t} =
  \sum_{i=1}^{N} \sum_{k=1}^{K} \sum_{j=1}^{K}
  &v_{kj}(\fatx_{i \cdot}+ \fatm_{kj})
  p(\fatx_{1 \cdot},\ldots, \fatx_{i \cdot}+\fatm_{kj},\ldots,
  \fatx_{N \cdot}, t)                                                     \\
  \label{eq:defdiff2} 
  -&v_{kj}(\fatx_{i \cdot})p(\fatx, t).
\end{align}
The corresponding transition vector $\fatm_{kj}$ is zero except for
two components: $m_{kj,k} = 1$ and $m_{kj,j} = -1$.

By combining (\ref{eq:defCME}) and (\ref{eq:defdiff2}), we arrive at
the RDME for a chemical system with reactions and diffusion
\begin{equation}
  \label{eq:defRDME}
  \begin{array}{llll}
  \displaystyle{\frac{\partial p(\fatx, t)}{\partial t} =
  \calM p(\fatx, t)+\calD p(\fatx, t)}.
  \end{array}
\end{equation}

We now give a few comments on the validity of the RDME. Denote the
molecular reaction radius by $\rho_R$ and the shortest average life
time of the molecular species by $\tau_{\min}$ \cite{BeHi85}. Then the
requirement for the size $h$ of a cell is
\begin{equation}
  \label{eq:hbound}
  \rho_R^2 \ll h^2 \ll \alpha\gamma \tau_{\min},
\end{equation}
where $\alpha$ is $\ordo{1}$ and depends on the cell shape and the
dimension \cite{ElE, VanKampen}. Firstly, the upper bound guarantees
that the mixing in a cell by diffusion is sufficiently fast for the
molecules to be homogeneously distributed there. Another
interpretation is that with slow diffusion, a better spatial
resolution is necessary since the solution becomes less
smooth. Secondly, there is also a \emph{lower} bound on $h$ for the
modeling at a meso level to be meaningful. This lower bound guarantees
that association and disassociation events can be properly localized
within a computational cell. The model breaks down if we let $h \to 0$
despite the fact that there is often a meaningful stochastic partial
differential equation (SPDE) in the limit. This SPDE, however, only
remains valid as the first term in a system size expansion
\cite[Ch.~8.2]{Gardiner}. It is understood with such an expansion that
we are working on a scale where the cell size appears small but is
still large enough to contain many molecules.

Let $\tri_k$ be a cell with a part of its boundary on
$\partial\Omega$. Since there is no transport of molecules out from
$\Omega$ at $\tri_k$ with the reactions in (\ref{eq:exdiff}), the
corresponding boundary condition at the macro level is a Neumann
condition. Most results for the numerical solution of diffusive
systems are derived under Dirichlet boundary conditions and we
therefore now discuss the associated conditions for the RDME. Suppose
that $\Omega$ is surrounded by a reservoir such that the number of
molecules is fixed at the boundary $\partial \Omega$. Then
$\fatx_{\cdot k}$ in those cells is given by the reservoir data. If
$x_{ij} = 0$, then $v_{jk} = 0$ in (\ref{eq:exdiff}) and there is no
diffusion of molecules from the boundary to the interior cells,
(cf. (\ref{eq:exdiff}))
\begin{equation}
  \label{eq:exdiff2}
  X_{ik}\rightarrow X_{ij} = \emptyset,\;
  v_{kj}(\fatx_{i \cdot}) = q_{kj} x_{ik}.
\end{equation}
The transition vector $\fatm_{kj}$ is zero except for $m_{kj,k} = 1$.

\subsection{Relation to macroscopic equations}

Define the concentrations $\phi_{ij}$ of species $i$ in cell
$\tri_{j}$ at a macroscopic level as the expected values of
$|\tri_j|^{-1}x_{ij}$. The RRE for $\phi_{ij}$ is obtained from the
CME (\ref{eq:defCME}) by deriving equations for the mean values of
$x_{ij}$ as in \cite{Gardiner, VanKampen}. The system of ODEs defining
the RRE is
\begin{align}
  \label{eq:defRRE}
  \frac{d \phi_{ij}}{dt} = \omega_{i}(\fatphi_{\cdot j}) \equiv 
  -\sum_{r=1}^{R} \frac{n_{ri}}{|\tri_{j}|}
  w_r(|\tri_{j}|\fatphi_{\cdot j}).
\end{align}

If there are no reactions but only diffusion, then a similar set of
macroscopic equations may be derived. From the similarity between
(\ref{eq:defCME}) and (\ref{eq:defdiff2}) we have
\begin{align}
  \nonumber
  \frac{d \phi_{ij}}{dt} &=
  -\sum_{k=1}^{K} \frac{m_{kj,j}}{|\tri_{j}|}
  v_{kj}(|\tri_{k}|\fatphi_{i \cdot})+
  \frac{m_{jk,j}}{|\tri_{j}|} v_{jk}(|\tri_{j}|\fatphi_{i \cdot})         \\
  \label{eq:defdiffE}
  &= \sum_{k=1}^{K} \frac{|\tri_{k}|}{|\tri_{j}|}q_{kj} \phi_{ik}-
  \left(\sum_{k=1}^{K} q_{jk}\right) \phi_{ij}.
\end{align}
The diffusion equation (\ref{eq:defdiffE}) has the form
\begin{equation}
  \label{eq:defdiffE2}
  \displaystyle{\frac{d \fatphi^T_{i \cdot}}{dt} =
  \gamma D \fatphi^T_{i \cdot}}
\end{equation}
for each $i$ (cf.~(\ref{eq:diffusion_comps})). The diffusion matrix
has the elements $\gamma D_{jk} = q_{kj}|\tri_{k}|/|\tri_{j}|,\; j
\neq k$, and $\gamma D_{jj} = -\sum_{k \not = j}q_{jk}$. 

The positive coefficients $q_{kj}$ in our diffusion model are non-zero
only if cell $k$ and $j$ share a common point in 1D, a common edge in
2D, or a common facet in 3D. The diffusion matrix $D$ therefore has a
sparsity pattern matching the connectivity of the partitioning of
$\Omega$ into computational cells.

Assume for now that the diffusion is isotropic on a Cartesian lattice
in 1D with constant cell size $h$ so that the probability $q_{kj}$ to
move from cell $\tri_{k}$ to cell $\tri_{j}$ is equal to the
probability $q_{jk}$ to move in the opposite direction. With isotropy,
$D$ has a particularly simple structure and is symmetric. For example,
if $\Omega = [0,1]$ with $h = 1/K$ and diffusion $q$ we have that
\begin{equation}
\label{eq:diff1D}
\begin{array}{llll}
  \dot{\phi}_{i1} &= q(-{\phi}_{i1}+{\phi}_{i2}), \\
  \dot{\phi}_{ij} &= q({\phi}_{i,j-1}-2{\phi}_{ij}+{\phi}_{i,j+1}),\; 
            j=2,\ldots,K-1, \\
  \dot{\phi}_{iK} &= q({\phi}_{i,K-1}-{\phi}_{iK}).
\end{array}
\end{equation}
Let $q = \gamma/h^2$ as in~(\ref{eq:diffusion_comps}) and let $h
\to 0$. Then the solution of (\ref{eq:diff1D}) converges to
the solution $\phi_i(x, t)$ of the diffusion equation in 1D
\begin{equation}
\label{eq:diffeq1D}
\begin{array}{llll}
   \displaystyle{\frac{\partial \phi_i}{\partial t}=
        \gamma \frac{\partial^2 \phi_i}{\partial x^2}},\quad
   \displaystyle{\frac{\partial \phi_i}{\partial x}(0, t)=
         \frac{\partial \phi_i}{\partial x}(1, t)=0}.
\end{array}
\end{equation}
If boundary values $\phi_{i1}$ and $\phi_{iK}$ are given, then the
solution converges to the solution of the PDE 
in (\ref{eq:diffeq1D}) with boundary conditions
$\phi_i(0, t) = \phi_{i1} = g_{i0}$ and $\phi_i(1, t) = \phi_{iK} =
g_{i1}$ for some $g_{i0}, g_{i1}\ge 0.$

Suppose that the interior of $\Omega = [0,1]\times [0, 1]$ is covered
by square cells of size $h\times h$ in 2D or that $\Omega=[0, 1]\times
[0, 1]\times [0, 1]$ in 3D and partitioned into cubic cells of size
$h\times h\times h$. In both cases, $D$ is symmetric and will
approximate the Laplacian $\Delta$ with Neumann boundary conditions.
With the normal derivative of $\phi_i$ at the boundary
$\partial\Omega$ written as $\partial \phi_i/\partial n$, the solution
$\phi_i$ converges when $h \to 0$ to the solution of
\begin{equation}
\label{eq:diffeq2D}
\begin{array}{llll}
   \displaystyle{\frac{\partial \phi_i}{\partial t}=
        \gamma  \Delta\phi_i}\;\; {\rm in}\;\; {\Omega},\;\;
   \displaystyle{\frac{\partial \phi_i}{\partial n}=0
     \;\; {\rm on}\;\; {\partial\Omega}},
\end{array}
\end{equation}
If $\phi_i$ is given data in the boundary cells, then the boundary
conditions in (\ref{eq:diffeq2D}) will be Dirichlet type, $\phi_i=g_i$
on $\partial \Omega$.

The macroscopic approximation of the diffusive part of
(\ref{eq:defRDME}) with a vanishing cell size in a Cartesian mesh thus
satisfies a diffusion equation (\ref{eq:diffeq1D}). The macroscopic
counterpart to $\calM$ in (\ref{eq:defRDME}) is the RRE
(\ref{eq:defRRE}). A macroscopic concentration $\phi_i$ affected by
both chemical reactions and diffusion fulfills a RDE
\begin{equation}
\label{eq:defRDE}
\displaystyle{\frac{\partial \phi_{i}}{\partial t} =
  \omega_i(\fatphi)+\gamma  \Delta\phi_i},\; i=1,\ldots,N.
\end{equation}

A Cartesian mesh for discretization of $\Delta \phi$ has many
advantages but is impractical for curved inner and outer boundaries of
$\Omega$. A RDME on unstructured meshes is proposed in 
Section~\ref{sec:diff}.

\subsection{Relation to microscopic equations}
\label{sec:FEM1D}

If we consider pure diffusion in 1D, it is possible to directly
compare the coefficients obtained from the finite element
discretization with the expected value of the first exit time of
Brownian motion from a finite interval. Using linear basis functions
in a FEM discretization of (\ref{eq:diffeq1D}) on a mesh with vertices
$x_j$ and cell sizes $h_j=x_j-x_{j-1}$, in the interior of $\Omega$
the non-zero entries of the stiffness matrix $S$ and the mass matrix
$M$ are

\begin{equation}\label{eq:1DSM}
\begin{array}{ccc}
\displaystyle
\gamma S_{j,j-1}=\frac{\gamma}{h_j},&
\displaystyle
\gamma S_{j,j+1}=\frac{\gamma}{h_{j+1}},&
\displaystyle
\gamma S_{jj}=-\frac{\gamma}{h_j} - \frac{\gamma}{h_{j+1}},\\
\displaystyle
M_{j,j-1}=\frac{1}{6}h_j,&
\displaystyle
M_{j,j+1}=\frac{1}{6}h_{j+1},&
\displaystyle
M_{jj}=\frac{1}{3}(h_j + h_{j+1}).
\end{array}
\end{equation}
\noindent
After mass lumping the coefficients corresponding to jumps from node
$j$ to its neighbors in \eqref{eq:diffusion_comps} are given by

\begin{equation}\label{eq:femjump}
\begin{array}{cc}
\displaystyle
q_{j,j-1}=\frac{2\gamma}{h_j(h_j+h_{j+1})},&
\displaystyle
q_{j,j+1}=\frac{2\gamma}{h_{j+1}(h_j+h_{j+1})}.
\end{array}
\end{equation}
\noindent
Consequently, the exponentially distributed waiting time for the next
event at node $j$ has the expected value
\begin{equation}\label{eq:exfem}
q_{j,j-1}+q_{j,j+1} = \frac{2\gamma}{h_jh_{j+1}}.
\end{equation}
On a uniform grid we recover the jump coefficients $\gamma/h^2$ and
the parameter $2\gamma/h^2$ used in \cite{HFE}.

Interpreted as the average time for a Brownian particle starting at
node $x_j$ to reach either of its neighbors, we can compare the value
in \eqref{eq:exfem} with the actual expected value of the first exit
time $\tau = \inf\{t: \xi_t \notin (x_{j-1},x_{j+1})\}$ where
$\xi_{t}$ is defined in (\ref{eq:Ito}). A straightforward application
of Dynkin's formula \cite[Ch.~7.4]{OKS} yields
\begin{equation}
  E_{x_j}[\tau]^{-1}=\frac{2\gamma}{h_j h_{j+1}},
\end{equation}
in accordance to \eqref{eq:exfem}. The probabilities to exit at
$x_{j-1}$ and $x_{j+1}$ respectively are given by
$h_{j+1}/(h_j+h_{j+1})$ and $h_{j}/(h_j+h_{j+1})$, and using this we
recover the jump coefficients \eqref{eq:femjump}. In this sense, the
coefficients of the mesoscale model obtained from the dicretization of
the macroscale equation is consistent with the microscale
description. It is worth noticing however, that $\tau$ is not
generally exponentially distributed \cite[p.~212]{BM}.

%% file: diff.tex
\section{Diffusion coefficients}
\label{sec:diff}

\begin{figure}[htp]
  \centering
  \includegraphics[width = 5.5cm]{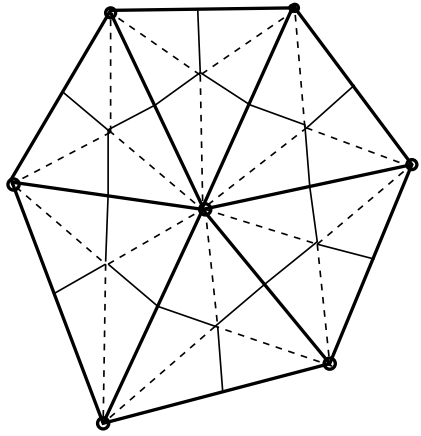}
  \caption{The primal mesh (thick lines) with the vertices (small circles),
  the dual mesh (thin lines) and the bisectors of the
  triangles (dashed lines).}
  \label{fig:unstruct}
\end{figure}

Consider a part of an unstructured mesh in 2D covering $\Omega$ with a
polygonal boundary $\partial\Omega$ in Figure~\ref{fig:unstruct}.  The
primal mesh consists of triangles with the vertices in the corners.
The cells $\tri_k$ in the dual mesh are polygons and in the interior
of $\Omega$, the center of $\tri_k$ is the vertex $k$. The edges of
the polygon coincide with a part of the bisectors of the triangles or
with the boundary $\partial\Omega$. The corners of an inner $\tri_k$
are the barycenters of the triangles and the midpoints of the edges
from its center vertex. In 1D, the cell in the primal mesh is a line
segment with a vertex in both ends and the dual mesh also consists of
line segments shifted with respect to the primal mesh and a vertex in
the center. The primal cell is a tetrahedron and the dual cell is a
polyhedron in 3D.

Similarly, the dual cells $\tri_k$ in a Cartesian structured mesh in
2D with vertices $(ih, jh), i=0,\ldots,K,\; j=0,\ldots,K,$ are the
cells in the staggered mesh. In the interior of $\Omega$, the vertices
are in the center of $\tri_k$.

A finite element discretization of the RDE in (\ref{eq:defRDE}) with
continuous, piecewise linear basis and test functions on the primal
unstructured mesh generates a system of equations to solve for the
nodal values $\fatphi$. The components of $\fatphi$ are
$\phi_{ij}(t)$, the concentration of species $i$ in vertex $j$ at time
$t$ \cite[Ch.~15]{Tho} but can also be interpreted as the mean value
of the concentration of species $i$ in the dual cell $\tri_j$. The
system of equations is \cite[Ch.~14]{Tho}
\begin{equation}
\label{eq:FEMeq}
    \hat{M}\dot{\fatphi}=\fatf(\fatphi)+\gamma \hat{S}\fatphi.
\end{equation}
The mass matrix $\hat{M}$ and the stiffness matrix $\hat{S}$ are
symmetric, $\hat{M}$ is positive definite and $\hat{S}$ is negative
semi-definite with Neumann boundary conditions and negative definite
with Dirichlet conditions. The reaction terms are represented by the
nonlinear term $\fatf$. If $h$ is a suitable measure of the sizes of
the triangles in the mesh, then the error in the FEM solution of the
Dirichlet problem in the $L_2$-norm is $\ordo{h^2}$
\cite[Ch.~14]{Tho}.

Order the unknowns in $\fatphi$ according to a linear index so that
\begin{align}
  \fatphi = (\fatphi_{1 \cdot}, \cdots, \fatphi_{i \cdot}, \cdots,
  \fatphi_{N \cdot})^{T}.
\end{align}
Then $\hat{M}$ and $\hat{S}$ are block diagonal matrices where the
blocks ${M}$ and ${S}$ are identical, small mass and stiffness
matrices. The system (\ref{eq:FEMeq}) is simplified by introducing
mass lumping of $M$ and $\fatf$. Let $\hat{A}$ be a diagonal matrix
with diagonal blocks ${A}$ with ${A}_{jj} = \sum_{k=1}^K {M}_{jk}$ and
let similarly $\bar{\fatf}$ be the lumped version of $\fatf$. In 1D,
${A}_{jj}$ is the length of the dual cell with vertex $j$ in the
center, see (\ref{eq:1DSM}). Similarly in 2D, ${A}_{jj}$ 
is the area $|\tri_j|$ of the dual
cell in Figure~\ref{fig:unstruct} with vertex $j$ in the center
\cite[Ch. 15]{Tho} and in 3D, ${A}_{jj}$ is the volume of $\tri_j$
\cite{ThW}. The simplified system (\ref{eq:FEMeq}) is now
\begin{equation}
  \label{eq:FEMeq2}
  \dot{\fatphi} = \hat{\fatomega}(\fatphi)+\gamma \hat{D}\fatphi, \quad
  \hat{\fatomega} \equiv \hat{A}^{-1}\bar{\fatf}, \; \hat{D} \equiv \hat{A}^{-1}\hat{S},
\end{equation}
which is our approximation of (\ref{eq:defRDE}). The solution
$\fatphi$ of (\ref{eq:FEMeq2}) is also generally second order accurate
in 2D \cite{NiT}.

Let $D$ denote a block on the diagonal of $\hat{D}$. Then $D =
A^{-1}S$ and an off-diagonal component $D_{jk}$ is non-zero only if
two vertices $j$ and $k$ are connected by an edge. With Neumann
conditions, the diagonal blocks $D$ satisfy $D_{jj} = -\sum_{k \not =
j} D_{jk}$ so that $\sum_{k=1}^K D_{jk} = 0$ for every $j$. In other
words, the constant vector $\fate_{1} = (1, 1,
\ldots, 1)^T$ is in the null-space of ${D}$ and is the right
eigenvector with eigenvalue zero. The corresponding left eigenvector
$\fate_2$ has the diagonal elements of ${A}$ as components, $e_{2 j} =
A_{jj}$. Let
\[
    \bar{x}_{ij}=\sum_{\fatx}x_{ij}p(\fatx, t)=|\tri_j|\phi_{ij}
\]
denote the expected value of the number of molecules of species $i$ in
cell $j$. In a system without reactions, by (\ref{eq:FEMeq2})
\begin{equation}
  \label{eq:phicons}
  \displaystyle{0=\gamma \fate_{2}^T\hat{D}\fatphi_{i \cdot}^T=
  \frac{d}{dt}\fate_{2}^T{\fatphi}_{i \cdot}^T=
  \frac{d}{dt}\sum_{j=1}^K A_{jj}\phi_{ij}=
  \frac{d}{dt}\sum_{j=1}^K |\tri_j|\phi_{ij}
  =\frac{d}{dt}\sum_{j=1}^K \bar{x}_{ij}}, 
\end{equation}
i.e.~the total number of molecules of each species are conserved by
the diffusion. This is not the case with Dirichlet conditions, where
$D_{jj} \le -\sum_{k \not = j} D_{jk}$ and $D$ is non-singular.

\begin{figure}[htp]
  \centering
  \subfigure[]{\includegraphics[width = 5.5cm]{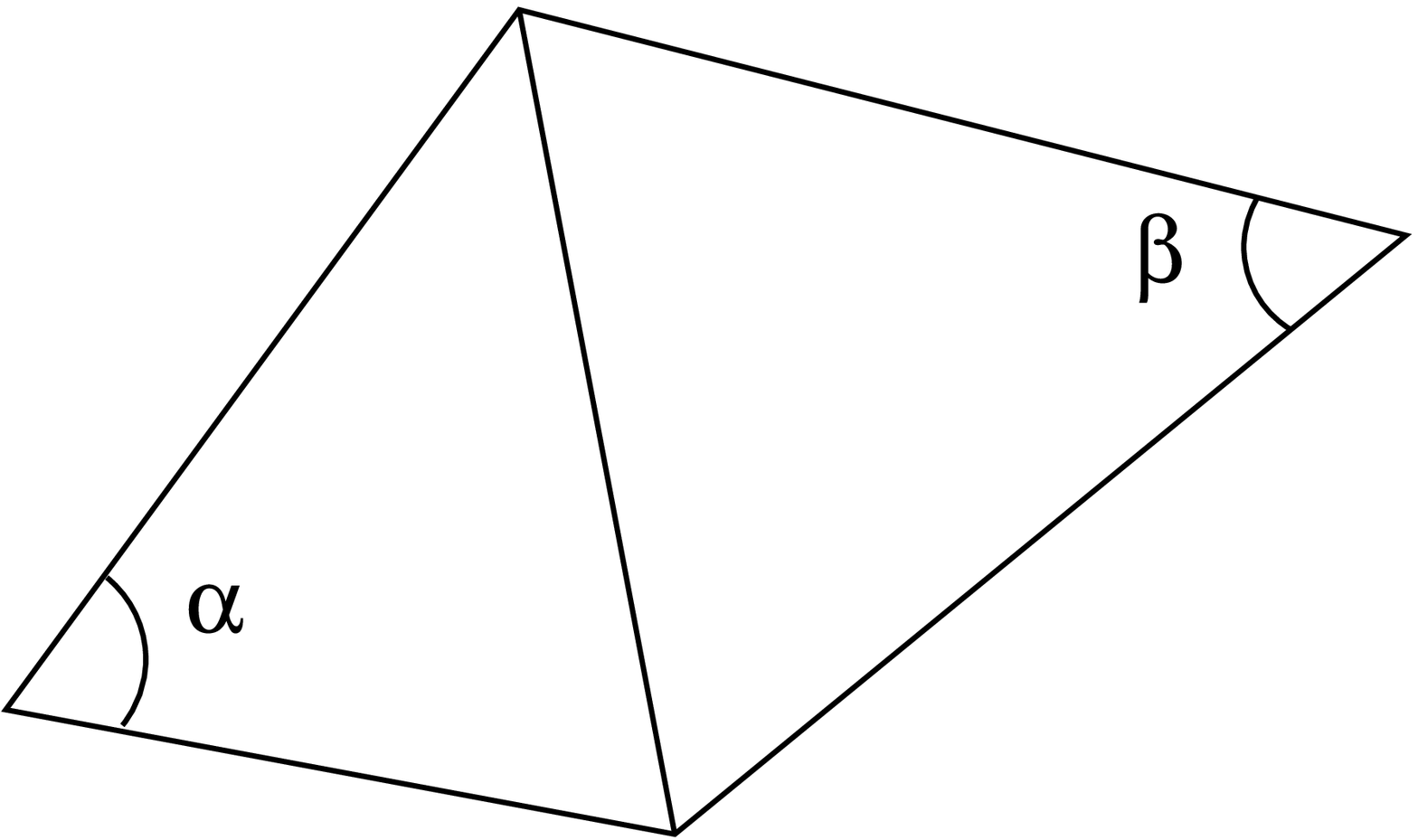}
        \hspace{10 mm}}
  \subfigure[]{\includegraphics[width = 4.5cm]{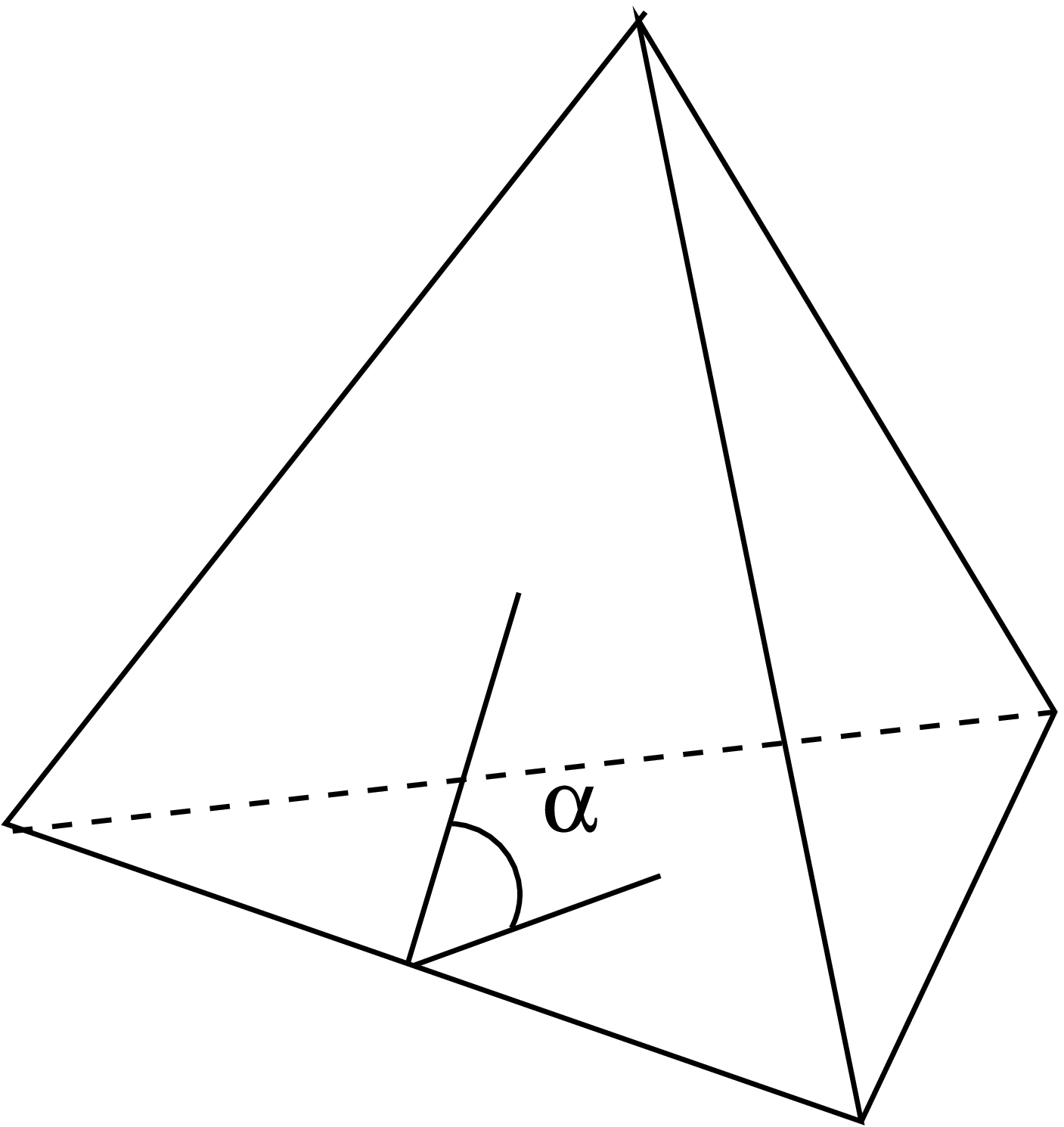}}
  \caption{The critical angles in 2D (a) and 3D (b) for positive
  off-diagonal elements in the stiffness matrix. }
  \label{fig:dihedral}
\end{figure}

By the calculation in Section~\ref{sec:FEM1D} in 1D, we find that
$D_{jk}>0$ if $k = j-1$ or $k = j+1$ and zero for the other
non-diagonal entries. If the mesh in 2D is a Delaunay triangulation
\cite{Geo}, then the Voronoi cells are close to the cells in the dual
mesh defined by Figure~\ref{fig:unstruct}. Two triangles share an edge
in Figure~\ref{fig:dihedral}a. The sum of the two opposing angles
$\alpha$ and $\beta$ in the triangles is less than or equal to $\pi$
in a Delaunay triangulation and $D_{jk}\ge 0$ when $j \neq k$
\cite{XuZ}.  In 3D, assume that the dihedral angle at an edge between
two facets of a tetrahedron, see Figure~\ref{fig:dihedral}b, is
non-obtuse $(\alpha\le \pi/2)$ for all tetrahedra in the mesh. Then
the elements of ${D}$ all have the right sign \cite{XuZ}. With these
properties of the discretization and without chemical reactions, a
discrete maximum principle is fulfilled
\cite{XuZ} for Dirichlet boundary conditions ensuring that with
non-negative initial conditions the solution remains non-negative. The
diffusion matrix ${D}$ generated by the FEM discretization of $\Delta$
in (\ref{eq:defRDE}) has the same properties as the diffusion matrix
in (\ref{eq:defdiffE2}). A review of these triangulations in three and
higher dimensions is found in \cite{BKKS}. They may not be trivial to
generate in higher dimensions than 2.

Discretization matrices of second order PDEs are frequently
\emph{$M$-matrices}~\cite[10.3]{Greenbaum}. Henceforth, we shall
assume the following slightly weaker property:
\begin{assumption}
\label{ass:M0}
The diffusion matrix $D$ for Dirichlet and Neumann boundary conditions
fulfills for $j \not = k$,
\[
  D_{jk} \ge 0,\; D_{jj} < 0,\; \sum_{k=1}^K D_{jk}\le 0.
\]
The last inequality is an equality for Neumann conditions.
\end{assumption}

The macroscopic elements in $\gamma {D}$ in (\ref{eq:FEMeq2}) define
the coefficients $q_{kj}$ in the mesoscopic model of diffusion
(\ref{eq:exdiff}). The diffusion matrix $D$ in (\ref{eq:defdiffE2})
has the form $\gamma D = A^{-1}Q A$, where $Q_{jk}=q_{kj}$ when $j\ne
k$. Since $D$ in (\ref{eq:defdiffE2}) and (\ref{eq:FEMeq2}) are
identical we have
\begin{equation}
\label{eq:Qdef}
    Q=\gamma S A^{-1}=\gamma D^T.
\end{equation}
Except for the case when all cells have the same size, $Q$ is 
generally unsymmetric.

The concentrations $\phi_{ij}$ defined as the expected values of
$x_{ij}/A_{jj}$ with the PDF in (\ref{eq:defRDME}) will satisfy
(\ref{eq:FEMeq2}). The molecules at the meso level can jump between
dual cells with a point (1D), edge (2D), or facet (3D) in common since
${D}_{jk}>0$ there. The connectivity graph of ${D}$ tells in which
direction a molecule can diffuse and the positive elements of ${D}$
are inversely proportional to the expected value of the first exit
time to leave the cell, cf.~(\ref{eq:diffusion_comps}) and
Section~\ref{sec:FEM1D}.

An alternative to the finite element discretization of the RDE in
(\ref{eq:defRDE}) is to use the finite volume method (FVM). Here, the
averages of the concentrations $\bar{\phi}_{ij}$ in the dual cells are
the degrees of freedom. Then
\begin{equation}
\label{eq:FVMeq}
\begin{array}{llll}
  \displaystyle{\frac{\partial \bar{\phi}_{ij}}{\partial t}}&=
     \displaystyle{\frac{1}{|{\tri}_{j}|}\int_{{\tri}_{j}} \frac{\partial {\phi}_i}{\partial t}\; dV
   =\frac{1}{|{\tri}_{j}|}\int_{{\tri}_{j}} 
          \omega_i(\fatphi) +\gamma \nabla\cdot \nabla \phi_i\; dV}\\
   &=\displaystyle{\frac{1}{|{\tri}_{j}|}\int_{{\tri}_{j}} \omega_i(\fatphi) \; dV+
    \frac{1}{|{\tri}_{j}|}\int_{{\partial\tri}_{j}}\gamma \hat{\fatn} \cdot \nabla \phi_i  \; dS},
\end{array}
\end{equation}
where Gauss' theorem has been used and $\hat{\fatn}$ is the normal of
${\partial\tri}_{j}$. The reaction term in ${\tri}_{j}$ is
approximated by $\omega_i(\fatphi_{\cdot j})$. The gradient $\nabla
\phi_i$ is needed on the boundary of ${\tri}_{j}$ and different
approximations are possible. A simple one is to let $\hat{\fatn}\cdot
\nabla \phi_i\approx (\phi_{ik}-\phi_{ij})/h_{jk}$ where $h_{jk}$ is
the distance between vertex $j$ and a neighboring vertex $k$. The
resulting diffusion term is in $\tri_j$,
\[
    \displaystyle{\sum_k c_{jk}(\phi_{ik}-\phi_{ij})/h_{jk}},
\]
with summation over those vertices with an edge between $j$ and $k$
and where $c_{jk}$ equals the size of the edges or facets of
$\partial\tri_j$ adjacent to the same edge. Then the coefficients in
this discretization are interpreted as the diffusion coefficients in
(\ref{eq:defdiff2}) in the same manner as in the FEM case. The
difference between (\ref{eq:FVMeq}) and (\ref{eq:FEMeq2}) lies in the
approximation of the diffusion term.

Convergence of the FVM to the analytical solution is proved for
certain discretizations of the gradient in \cite{EGH, Her} but the
quality of the approximation seems to depend critically on the quality
of the mesh \cite{SGN}. This is one of the reasons why we prefer the
FEM approach. Another reason is that FVM is perhaps more suitable for
problems dominated by convection while chemical systems from molecular
biology tend to be of diffusive character.

%% file: mom.tex
\section{Moments of the diffusion}
\label{sec:mom}

The purpose of this section is to prove that the diffusion can be
accurately evolved deterministically when the number of diffusing
molecules is sufficiently large. Consider the equations for the
moments of $\fatx$ in a system with diffusion and Neumann boundary
conditions but without chemical reactions. The expected values and the
covariance matrices satisfy systems of ODEs. These equations are
derived in \cite{master_moment, VanKampen} for general propensities,
but with
\[
  v_{kj}(\fatx_{i \cdot})=Q_{jk} x_{ik}, \; j\ne k,
\]
being linear (cf.~(\ref{eq:exdiff})), they have a particularly simple
structure.

Applying the formulas in \cite{master_moment,FLH} or invoking
(\ref{eq:defdiffE2}), the first moment of the number of molecules of
species $i$ in a cell is given by
\begin{equation}
  \label{eq:avereq}
  \dot{\bar{\fatx}}_{i \cdot}^T=Q \bar{\fatx}_{i \cdot}^T.
\end{equation}
This equation is exact since the diffusion propensities $v_{kj}$ are
linear in $\fatx$ and no coupling to higher order moments exists.

The second moments or covariances of any species $i$ between cells $j$
and $k$ are denoted by $C_{jk}$. The equation for $C_{jk}$ is
\begin{equation}
\label{eq:coveq}
   \dot{C}_{jk}=\displaystyle{\sum_{l=1}^K Q_{kl}C_{jl} +
                             \sum_{l=1}^K Q_{jl}C_{kl}+F_{jk}},\; j,k=1,\ldots,K,
\end{equation}
with the driving term $F$ defined by
\[
  F_{jk}=\displaystyle{\sum_{\alpha=1}^K \sum_{\beta=1}^K 
         m_{\alpha \beta, j} m_{\alpha \beta, k} v_{\alpha \beta}(\bar{\fatx}_{i \cdot})}. 
\]
This can be written in matrix form since $C$ is symmetric,
\begin{equation}
\label{eq:covmateq}
   \dot{C}=CQ^T +QC+F.
\end{equation}
Using the properties of the diffusion propensities, the elements of
$F$ are
\begin{equation}
\label{eq:Wdef}
\begin{array}{llll}
  &\displaystyle{F_{jj}= \sum_{l=1, l\ne j}^K
  Q_{jl}\bar{x}_{il}}+Q_{lj}\bar{x}_{ij},\;
  &\displaystyle{F_{jk}= -(Q_{jk}\bar{x}_{ik}+Q_{kj}\bar{x}_{ij}),\; j\ne k},
\end{array}
\end{equation}
The covariance equation is also exact for the same reason as before.
 
The solution of (\ref{eq:avereq}) can be written
\begin{equation}
\label{eq:expecexct}
   \bar{\fatx}_{i \cdot}(t)=\bar{\fatx}^0_{i \cdot}\exp(Q^T t),
\end{equation}
where $\bar{\fatx}^0_{i \cdot}$ is the initial value at $t = 0$. The
eigenvalues of ${Q}$ and ${D}$ are all negative except for one which
is zero. The corresponding left eigenvector of $Q$ is $\fate_1$ and
the right eigenvector is $\fate_2$, both of them defined in
Section~\ref{sec:diff}. Hence,
\begin{equation}
\label{eq:expect}
   \bar{\fatx}_{i \cdot}(t)=\kappa_i\fate_2^T +\delta \bar{\fatx}_{i \cdot}(t).
\end{equation}
with an upper bound on $\|\delta \bar{\fatx}_{i \cdot}(t)\|$ given by
$c_\delta \exp(\lambda_2 t)$ where $\lambda_2$ is the negative
eigenvalue of $Q$ with smallest magnitude and the norm is the $\ell_2$-norm. Therefore,
\begin{equation}
\label{eq:limexpect}
   \lim_{t\rightarrow \infty}\bar{\fatx}_{i \cdot}(t)=\kappa_i\fate_2^T.
\end{equation}
By (\ref{eq:phicons}), $\bar{\fatx}_{i \cdot}\fate_1$ is constant and
we obtain
\begin{equation}
\label{eq:kappa}
  \kappa_i\fate_2^T\fate_1=\bar{\fatx}^0_{i \cdot}\fate_1 \Rightarrow
  \kappa_i=\displaystyle{\sum_{j=1}^K \bar{\fatx}^0_{ij}/\sum_{j=1}^K A_{jj}}.
\end{equation}

The explicit solution of (\ref{eq:covmateq}) is
\begin{equation}
\label{eq:covexct}
   C(t)=\exp(Q t)C_0\exp(Q^T t)+
        \int_0^t\exp(Q (t-s))F\exp(Q^T (t-s))\; ds,
\end{equation}
where $C_0$ is the initial value of the covariance at $t=0$.  Using
(\ref{eq:Wdef}), (\ref{eq:expect}), (\ref{eq:limexpect}), and
(\ref{eq:Qdef}) we find that when $t \to \infty$
\[
\begin{array}{llll}
  F_{jj}&=\gamma\kappa_i \sum_{l=1, l\ne j}^K
  A_{ll}S_{jl}/A_{ll}+A_{jj}S_{lj}/A_{jj}=
  2 \gamma\kappa_i\sum_{l=1, l\ne j}^K S_{jl}=-2\gamma\kappa_iS_{jj},\\
  F_{jk}&=-\gamma\kappa_i(A_{kk}S_{jk}/A_{kk}+A_{jj}S_{kj}/A_{jj})=
  -2\gamma\kappa_iS_{jk},\; j\ne k,\\
\end{array}
\]
and hence that,
\begin{equation}
\label{eq:limW}
   F = -2 \gamma\kappa_iS+\delta F,
\end{equation}
where $\|\delta F\|$ is bounded by $c_{\delta F} \exp(\lambda_2t)$. A
bound on the covariance matrix for finite time is derived in the next proposition.

\begin{proposition}
\label{prop:Cbound}
Suppose that $C_0=0$. Then 
\begin{equation}
  \label{eq:Wbnd}
  \displaystyle{\|C(t)\|\le
  c_{F}\frac{\max_j{A_{jj}}}{\min_j{A_{jj}}}
  \int_0^t \|\bar{\fatx}_{i \cdot}\|\; ds},
\end{equation}
for some bounding constant $c_{F}$ such that $\|F\| \le c_{F}
\|\bar{\fatx}_{i \cdot}\|$.
\end{proposition} 

\begin{proof}
Let $\tilde{S} = A^{-1/2}SA^{-1/2}$. Then
\[
   \exp(Qt) = \exp(\gamma SA^{-1}t) = 
  A^{1/2}\exp(\gamma \tilde{S} t)A^{-1/2}.
\]
The symmetric matrix $\gamma\tilde{S}$ has the same eigenvalues as $Q$. Let
the unitary matrix $U$ have the eigenvectors of $\tilde{S}$ as columns and let $\Lambda$
have the corresponding eigenvalues on the diagonal. A bound on the
exponential of $Q\tau$ with $\tau\ge 0$ is
\begin{align*}
  \|\exp(Q\tau)\| &\le \|A^{1/2}\| \|A^{-1/2}\|
  \|U\exp(\gamma\Lambda \tau)U^T\|                                      \\
  &\le
  \|\exp(\gamma\Lambda \tau)\| \max_j\sqrt{A_{jj}}/\min_j\sqrt{A_{jj}}
  \le \max_j\sqrt{A_{jj}}/\min_j\sqrt{A_{jj}}.
\end{align*}
The same bound is valid for $\exp(Q^T\tau)$. Hence,
\begin{align*}
  &\left\| \int_0^t \exp(Q (t-s))F\exp(Q^T (t-s)) \; ds \right\|        \\
  &\le \int_0^t \|F\|\max_j{A_{jj}}/\min_j A_{jj}\; ds \le 
  c_{F} \frac{\max_j{A_{jj}}}{\min_j A_{jj}}
  \int_0^t \|\fatx_{i \cdot}\|\; ds.
\end{align*}
\end{proof}

\begin{remark}
Using (\ref{eq:covexct}) and (\ref{eq:limW}), one can show that
$\|C(t)\|$ is bounded when $t\rightarrow \infty$. If $C_0=0$, then
$\|C(t)\|\sim \kappa_i$ for large $t$.
\end{remark}

It follows from the proposition that the variance of $\fatx_{i \cdot}$
is proportional to $\|\bar{\fatx}_{i\cdot}\|$ in a bounded time
interval so that the standard deviation is proportional to
$\sqrt{\|\bar{\fatx}_{i \cdot}\|}$. When $t$ is large then
$\bar{\fatx}_{i \cdot}$ in (\ref{eq:limexpect}) is of the same order
as $\kappa_i$, a weighted average of the initial $\bar{\fatx}^0_{i
\cdot}$ (\ref{eq:kappa}). The standard deviation is proportional to
$\sqrt{\kappa_i}$ according to the remark after
Proposition~\ref{prop:Cbound}. Therefore, the quotient between the
standard deviation and the expected value is small for large copy
numbers implying that the expected value is indeed a good
approximation. Conversely, if the number of molecules in a cell is
small, then a description in terms of expectation values should not be
used.

%% file: hybrid.tex
\section{Time integration and hybrid diffusion}
\label{sec:hybrid}

For a discretization parameterized by the cell size $h$, the time to
compute a trajectory of a system with SSA spent in the diffusion part
of (\ref{eq:defRDME}) is proportional to $x \gamma/h^{2}$, where $x$
is the total number of diffusing molecules (cf.~(\ref{eq:exdiff}) and
(\ref{eq:diffusion_comps})). It follows that diffusion is the only
event with a \emph{total} intensity that increases with increasing
spatial resolution. In order to avoid this, we propose to split the
diffusion operator $\calD$ into two parts with the species with low
copy numbers in $\calD_L$ and the species with high copy numbers in
$\calD_H$. The diffusion in $\calD_H$ can then be advanced in time
macroscopically.

Order the species $X_i$ such that $X_i,\; i=1,\ldots,N_L,$ have low
copy numbers and $X_i,\; i=N_L+1,\ldots,N,$ have high numbers and let
\begin{align}
  \nonumber
  &\xi_i \equiv \displaystyle{\sum_{k=1}^K
  \sum_{j=1}^K}  \displaystyle{v_{kj}(\fatx_{i \cdot}+ \fatm_{kj})
  p(\fatx_{1 \cdot},\ldots, \fatx_{i \cdot}+\fatm_{kj},\ldots,
  \fatx_{N \cdot}, t)-v_{kj}(\fatx_{i \cdot})p(\fatx, t),}              \\
  \label{eq:defDD}
  &\calD_L p(\fatx, t) \equiv \displaystyle{\sum_{i=1}^{N_L}\xi_i},\;
  \calD_H \equiv \calD-\calD_L.
\end{align}
The operator on the right hand side of (\ref{eq:defRDME}) can be
written
\begin{equation}
\label{eq:defRDME2}
\begin{array}{llll}
\displaystyle{\frac{\partial p(\fatx, t)}{\partial t} =
  [\calM+\calD_L] p(\fatx, t)+\calD_H p(\fatx, t)}.
\end{array}
\end{equation}

It follows from Section~\ref{sec:mom} that the effect of diffusion on
the species with many molecules in each cell is well approximated by
the mean-field equations. In the numerical solution procedure, the
second part of (\ref{eq:defRDME2}) is therefore first advanced half a
time step $\Delta t/2$ with the macroscopic diffusion. Then the first
part is integrated a full step $\Delta t$ and finally the macroscopic
diffusion is applied for half a time step again. This is the Strang
splitting procedure \cite{Str} to solve (\ref{eq:defRDME2}) and below
we give conservative conditions under which both steps preserve the
non-negativity of the solution. 

For many relevant cases, a single trajectory gives sufficient insight
into the stochastic reaction-diffusion system but the same procedure
also works well for simultaneous simulation of $M$ trajectories. An
application could be to approximate the PDF by following an ensemble
of trajectories with state vectors $\fatx^m(t),\; m =
1,\ldots,M$. Then $p$ is reconstructed according to
\begin{equation}
\label{eq:papprox}
   p(\fatx, t^n)\approx \frac{1}{M}\sum_{m=1}^M \Psi^m,\qquad 
   \Psi^m=\left\{\begin{array}{l}1, \;\fatx^m(t^n)=\fatx ,\\
                                 0, \;{\rm otherwise}.\end{array}\right.
\end{equation}

In order to advance the trajectories $\Delta t$ in time from $t^n$ to
$t^{n+1}$, each trajectory is first integrated in time from $t^n$ to
$t^n+\Delta t/2$ by solving (\ref{eq:avereq}) 
for the species $i=N_L+1,\ldots,N$. The time derivative in
(\ref{eq:avereq}) is discretized by the trapezoidal (or
Crank-Nicolson) method of second order temporal accuracy and the new
state $\fatx_{i \cdot}^{n+1/4}$ for each trajectory is the solution of
\begin{equation}\label{eq:trap}
   \displaystyle{(I-\frac{1}{2}\Delta t Q)(\fatx_{i \cdot}^{n+1/4})^T =
   (I+\frac{1}{2}\Delta t Q)(\fatx_{i \cdot}(t^{n}))^T},\; i=N_L+1,\ldots,N.
\end{equation}
Alternatively, a scheme with an error $\ordo{\Delta t}$ is the Euler
backward method
\begin{equation}\label{eq:Eulerbw}
   \displaystyle{(I-\frac{1}{2}\Delta t Q)(\fatx_{i \cdot}^{n+1/4})^T=
   (\fatx_{i \cdot}(t^{n}))^T}.
\end{equation}

The time step in the Strang splitting must be sufficiently small to
resolve the shortest time scale of the reactions $\tau_{\min}$ in
(\ref{eq:hbound}). Thus, we can take $\Delta t\sim \tau_{\min}$ and by
(\ref{eq:hbound}) and \cite{ElE}
\begin{align}
  \label{eq:timstlow}
  h^{2} &\ll \alpha \gamma \Delta t.
\end{align} 
For such a time step, an explicit method for integration of
(\ref{eq:avereq}) would be very inefficient. Only an implicit method
suitable for stiff problems, such as the trapezoidal method, can
efficiently advance the solution of (\ref{eq:avereq}) in time in a
stable manner.

The next step is to evolve the $M$ trajectories with SSA
\cite{gillespie} for one time step $\Delta t$ using the reduced master
equation
\begin{equation}\label{eq:reducME}
\displaystyle{\frac{\partial p(\fatx, t)}{\partial t} = 
   \calM p(\fatx, t)+\calD_L p(\fatx, t)}.
\end{equation}
If $\Delta t$ is short then there may be no events in many of the
realizations of the process and some computational work will be
wasted.
  
The final step is to evolve each trajectory $\fatx^{m}_{i \cdot}$ half
a time step again using the macroscopic diffusion. The result is an
approximative sample of $p$ in (\ref{eq:defRDME2}) at time $t^{n+1}$.

The error due to the Strang splitting and the time discretization by
the trapezoidal method is $\ordo{\Delta t^{2}}$. Similarly, the error
due to the FEM approximation of the diffusion is $\ordo{h^2}$
\cite[Ch. 7]{Tho}.
 
It is more difficult to estimate the error induced from using the
macroscopic diffusion without being overly pessimistic. A bound on the
local \emph{single trajectory} error can be obtained as follows. Write
$\|C(\Delta t)\| = \ordo{\Delta t \|\bar{\fatx}_{i \cdot}\|}$ by
Proposition~\ref{prop:Cbound}. Then the local stochastic error
(standard deviation) relative to $\bar{\fatx}_{i \cdot}$ is
$\ordo{\Delta t^{1/2} \|\bar{\fatx}_{i \cdot}\|^{-1/2}}$ for $i >
N_{L}$ since only species with high copy numbers participate in the
macroscopic diffusion. This simple bound, however, gives no global
estimate since reasonable and sufficient stability properties of the
system are difficult to prescribe. When using deterministic diffusion
for some species it must simply be regarded as given \textit{a priori}
that the diffusion noise for those variables has little or no impact
on the system as a whole. Note that computing averages in order to
approximate expected values will generally make this error
substantially smaller since the macroscopic diffusion is
\emph{exact} in expectation (cf.~(\ref{eq:avereq})).

In a system with only diffusion, the algorithm has the following two
properties.

\begin{proposition}
\label{prop:nonneg}
Assume that $D$ satisfies Assumption~\ref{ass:M0}, that $\fatx^m(0) \ge
0$ for all trajectories $m = 1,\ldots,M$, and that $\Delta t
\le h_{\min}^2/6\gamma$ for the trapezoidal method in (\ref{eq:trap}),
where $h_{\min}$ is the minimal distance between a vertex and the
opposing edge in a triangle in the mesh. Then in a system without
chemical reactions, the copy numbers in the trajectories computed by
the hybrid algorithm remain non-negative for $t > 0$. For the Euler
backward method (\ref{eq:Eulerbw}), there is no time step restriction
for non-negativity.
\end{proposition}

\begin{proof}
If ${D}$ satisfies Assumption~\ref{ass:M0}, then $Q$ satisfies the
assumption by (\ref{eq:Qdef}). Let $({\fatx}_{i
\cdot}(t^n))^T=A^{1/2}\faty^n$ in a symmetrization of
(\ref{eq:trap}). Then
\[
     \displaystyle{(I- \tilde{S})
     \faty^{n+1/4}=
   (I+\tilde{S})\faty^{n}=\fatg},\; 
    \tilde{S}=0.5\gamma\Delta t A^{-1/2}SA^{-1/2}.
\]
The symmetric matrix $\tilde{S}$ has components $\tilde{S}_{jk}\ge
0,\; j\ne k,$ and $\tilde{S}_{jj}< 0$ and $I-\tilde{S}$ is positive
definite since $\tilde{S}$ is negative semi-definite. By \cite[Lemma
15.4]{Tho}, $(I-\tilde{S})^{-1}_{jk}\ge 0$ and by \cite[Theorem
15.6]{Tho} if $\Delta t\le h_{\min}^2/6\gamma$ and $\faty^{n}\ge 0$,
then the right hand side $\fatg$ is non-negative. Consequently,
\[
   \faty^{n+1/4}=(I-\tilde{S})^{-1}\fatg\ge 0,
\]
and therefore $\bar{\fatx}_{i \cdot}^{n+1/4}\ge 0$. The only
differences for the Euler backward method are that the right hand side
$\fatg$ equals $\faty^{n}$ and is always non-negative and $\tilde{S}$
is twice as large.
The intermediate SSA-step also preserves the non-negativity as do the
final step. Thus, the copy numbers of the species computed by the
hybrid algorithm remain non-negative.
\end{proof}

\begin{remark}
The upper bound on $\Delta t$ in the proposition is quite restrictive
considering the requirements for resolving the reactive time scale in
the splitting in (\ref{eq:timstlow}). In practice the solutions in the
macroscopic diffusion step stay non-negative with much longer $\Delta
t$ since the mean values are large for the species involved in
$\calD_H p$.
\end{remark}

\begin{proposition}
\label{prop:cons}
In a system with only diffusion and $\sum_k D_{jk}=0,\; j=1,\ldots,K$,
the total number of molecules of each species in a trajectory is
constant.
\end{proposition}

Note that the exact solution to the equations for the concentrations
has the same property in (\ref{eq:phicons}).

\vspace{3mm}

\begin{proof}
The vector $\fate_1$ satisfies $D\fate_1=0$ and $\fate_1^TQ =
\fate_1^TD^T = 0$. In the first step of the hybrid algorithm
(\ref{eq:trap}), we have
\[
   \displaystyle{\fate_1^T(I-\frac{1}{2}\Delta t Q)({\fatx}_{i \cdot}^{n+1/4})^T=
   \fate_1^T({\fatx}_{i \cdot}^{n+1/4})^T=
   \fate_1^T(I+\frac{1}{2}\Delta t Q)({\fatx}_{i \cdot}(t^{n}))^T=
   \fate_1^T({\fatx}_{i \cdot}(t^{n}))^T}.
\]   
Consequently, 
\[
   s_i^{n+1/4}\equiv\sum_{m=1}^M\sum_{j=1}^K x_{ij}^{m,n+1/4}=
   \sum_{m=1}^M\sum_{j=1}^K x_{ij}^{m,n}\equiv s_i^{n},\; i=N_L+1,\ldots,N,
\]
and $s_i$, the sum of the copy number of species $i$ over all cells, 
is thus preserved by (\ref{eq:trap}). This sum over
all cells is also preserved by diffusion simulated by SSA in every
trajectory in the intermediate step in the hybrid algorithm. Finally,
$s_i$ is preserved in the last step of the Strang splitting.
\end{proof}

The accuracy and efficiency of the algorithm are evaluated in the next
section.

%% file: res.tex
\section{Numerical results}
\label{sec:res}

The algorithm for the RDME in Section~\ref{sec:hybrid} is applied to
the diffusion equation and to two different systems from molecular
biology. The convergence of samples from the mesoscopic diffusion
model to the solution of the macroscopic equation is illustrated in
Section~\ref{sec:diffres}. The method is applied to a model of a
bi-stable reaction network in Section~\ref{sec:ehrelf}. Finally, in
Section~\ref{sec:hybridres} we illustrate the potential of the hybrid
method by comparing it to a purely stochastic simulation. The meshes,
the stiffness and the mass matrices are generated using the
PDE-toolbox~\cite{PDEtool} in MATLAB.

\subsection{Diffusion}
\label{sec:diffres}

\begin{figure}[htp]
  \centering
  \includegraphics[width = 4.5cm]{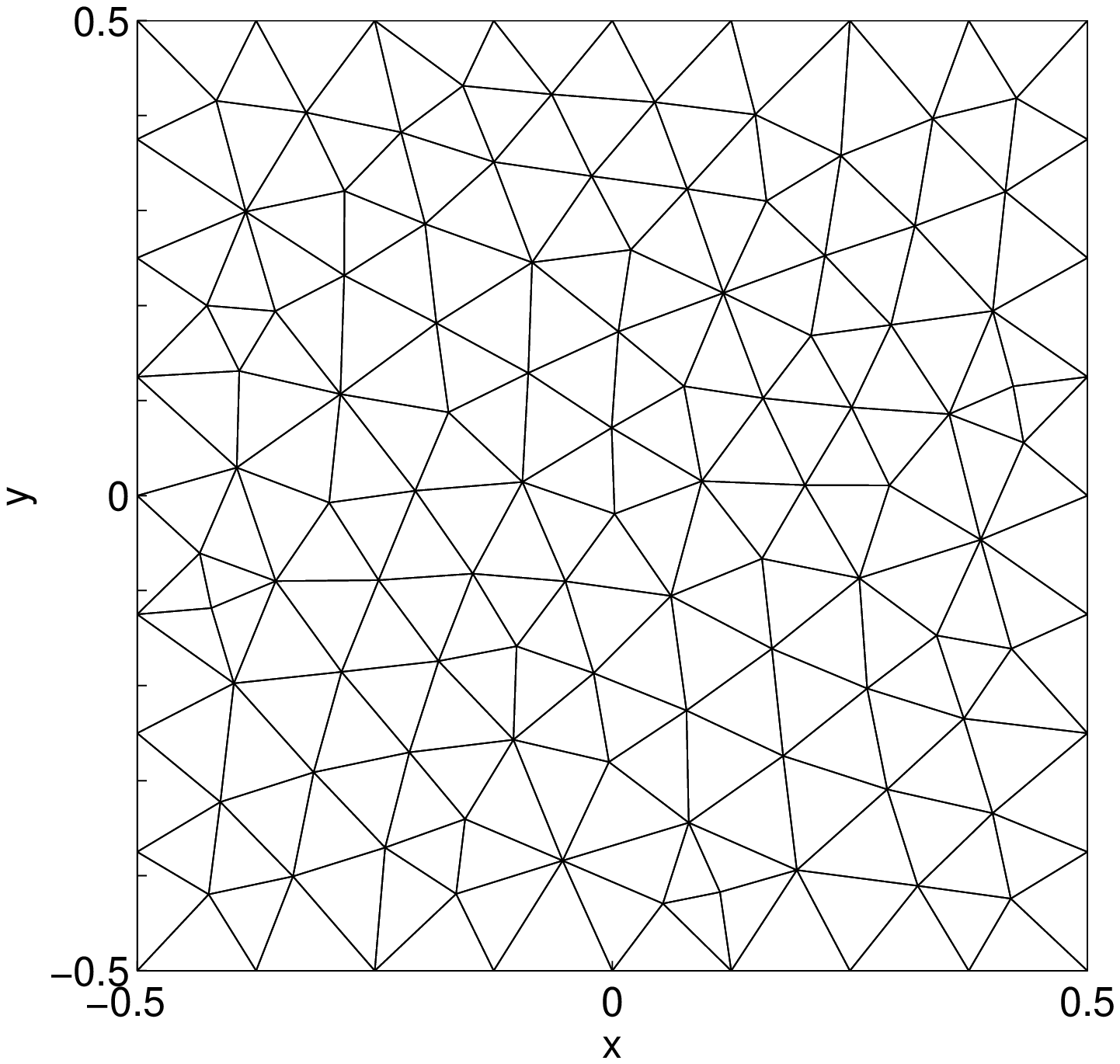}\hspace{5 mm}
  \includegraphics[width = 5.5cm]{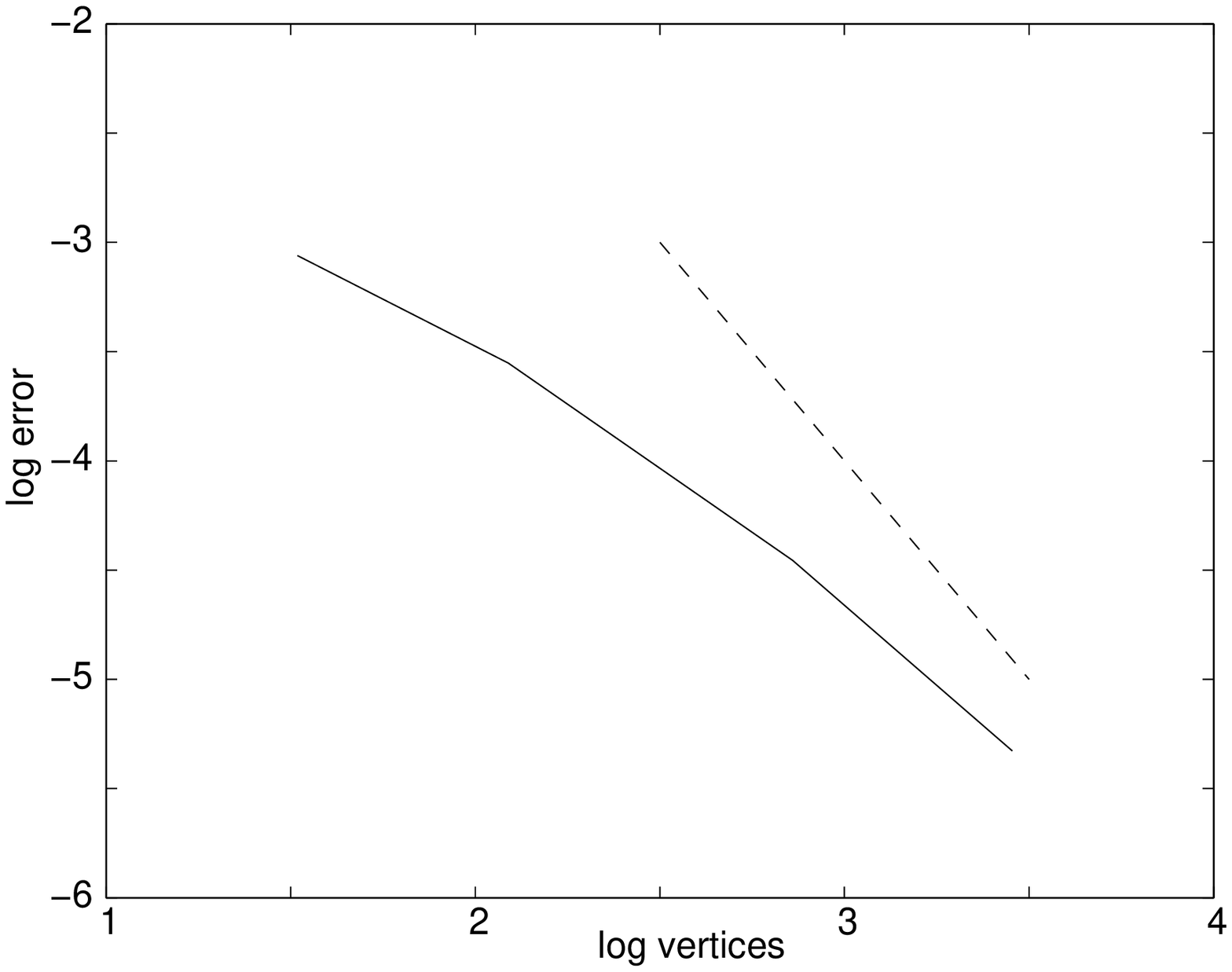}
  \caption{The triangular mesh with 123 vertices (left) and the logarithm 
  of the error vs.~the number of vertices (right). The $\ell_2$-norm of the
  error (solid) is compared to the asymptotic rate of convergence
  (dashed).}
  \label{fig:unstruct2}
\end{figure}

The diffusion equation with Neumann boundary conditions and initial data
\[
\begin{array}{lll}
   u_t=\gamma(u_{xx}+u_{yy}),\; {\rm in}\; \Omega=[-0.5, 0.5]\times[-0.5, 0.5],\; 
\gamma=10^{-3},\\
   \displaystyle{\frac{\partial u}{\partial n}}=0,\; {\rm on}\; \partial\Omega,\;
   u(x, y, 0)=100(1-\cos(2\pi x)),
\end{array}
\]
has the analytical solution
\begin{equation}
  \label{eq:analyt}
  u_a(x, y, t)=100(1-\cos(2\pi x) \exp(-4\gamma\pi^2t)).
\end{equation}

The solution $u_d$ is computed with a FEM discretization of the space
derivatives with mass lumping and integrated in time by the
trapezoidal method as in (\ref{eq:trap}). The error $u_d-u_a$ in the
vertices in the $\ell_2$-norm is $\ordo{\Delta t^2}+\ordo{h^2}$, see
Section~\ref{sec:hybrid}. The behavior of the spatial error at $t = 1$
is confirmed in Figure~\ref{fig:unstruct2}. Asymptotically, when the
number of vertices increases the theoretical rate is obtained.

A stochastic simulation of the diffusion with $m$ molecules in the
mesh is compared with the analytical and the FEM solutions in
Table~\ref{tab:diff}. Two meshes are used: one with 33 vertices and a
maximum length $h_{max}$ of an edge of a triangle equal to 0.5 and one
with 123 vertices and $h_{max}=0.25$ (see
Figure~\ref{fig:unstruct2}). Each stochastic simulation starts with
100 molecules distributed according to $u(x, y, 0)$.  The average
concentrations $u_m$ in the cells $\tri_k$ are computed for $M$
trajectories such that $m=100M$. Then $u_m$ is compared to $u_a$ and
$u_d$ in the vertices at $t=1$. The weighted vector norms in $\ell_2$
and $\ell_\infty$ are defined by
\begin{equation}
  \label{eq:lnorms}
  \| u\|^2_2=\sum_j u_j^2 |\tri_j|,\; \| u\|_\infty=\max_j |u_j|.
\end{equation}
The differences in these norms divided by the system size 100, 
$\delta_a=(u_m-u_a)/100$ and $\delta_d=(u_m-u_d)/100$,  
for two different discretizations
are collected in Table~\ref{tab:diff}.

\begin{table}[htp]
\begin{center}
\begin{tabular}{|c|cc|cc|cc|cc|}\hline
&\multicolumn{4}{c|} {$\ell_2$} &
\multicolumn{4}{c|} {$\ell_\infty$}\\ \hline
&\multicolumn{2}{c|}{$h_{max}=0.5$} & \multicolumn{2}{c|}{$h_{max}=0.25$}&
\multicolumn{2}{c|}{$h_{max}=0.5$} & \multicolumn{2}{c|}{ $h_{max}=0.25 $}\\
  $\log_{10}{m}$ & $\delta_a$ & $\delta_d$ &$\delta_a$ & $\delta_d$ 
            &  $\delta_a$  &  $\delta_d$  &  $\delta_a$  &  $\delta_d$\\ \hline
 2&   .038    & .038   & .1     & .1     &   .52     & .53   &  4.5   & 4.5 \\
 3&   .019    & .018   & .029   & .029   &   .32     & .31   &  1.2   & 1.2  \\
 4&   .0053   & .0053  & .0073  & .0073  & .091    & .096  &  .23   & .23 \\
 5&   .0018   & .0016  & .0029  & .0030  &  .031    & .017  &  .16   & .17 \\
 6&   .0011   & .00057 & .0011  & .0010  & .016    & .0093 &  .047  & .042\\
 7&   .00083  & .00015 & .00039 & .00027 &  .013    & .0023 &  .014  & .010\\  \hline
\end{tabular}
\end{center}
\caption{The relative difference between the stochastic solution and 
the analytic solution $\delta_a$ or the FEM solution $\delta_d$ for different 
mesh sizes $h_{max}$ and total number of molecules $m$.}
\label{tab:diff}
\end{table} 

The error is expected to behave as $\ordo{h_{max}^2}+\ordo{m^{-1/2}}$
and this is what we observe in the table. The difference between $u_m$
and $u_d$ decays with increasing $m$ at the correct rate in both
norms. In the example in Figure~\ref{fig:unstruct2}, the
$\ell_2$-error is $8.7\cdot 10^{-4}$ when $h_{max} = 0.5$ and
$2.8\cdot 10^{-4}$ when $h_{max} = 0.25$ explaining the difference
between $\delta_a$ and $\delta_d$ in Table~\ref{tab:diff}.  When $m$
is large then the dominant term in $\delta_a$ is the discretization error.

\subsection{Domain separation in a bi-stable system}
\label{sec:ehrelf}

In this section we simulate a model of a bi-stable system, previously
investigated using the freely available software MesoRD~\cite{HFE}
in~\cite{ElE}. The model consists of eight chemical species
participating in twelve reactions, see Table \ref{tab:bistab}. Being
based on a double negative feedback mechanism, in the spatially
homogeneous case the system switches between states with mostly $A$
molecules and states where $B$ is dominating. The model is used to
illustrate and explain the observation that global bi-stability can be
lost in a spatially dependent system due to domain separation, when
the diffusion is slow.

\begin{table}[htp]
\begin{tabular}{llll}
  $E_A \xrightarrow{k_1} E_A + A$ & $E_B \xrightarrow{k_1} E_B + B $ &
  $E_A+B \overset{k_a}{\underset{k_d}{\rightleftharpoons}} E_AB$ &
  $E_B+A \overset{k_a}{\underset{k_d}{\rightleftharpoons}} E_BA$        \\
  $E_AB + B \overset{k_a}{\underset{k_d}{\rightleftharpoons}} E_AB_2$ &
  $E_BA \overset{k_a}{\underset{k_d}{\rightleftharpoons}} E_BA_2$ &
  $A \xrightarrow{k_4} \emptyset$ & $B \xrightarrow{k_4} \emptyset$
\end{tabular}
\caption{The chemical reactions of the bi-stable model. The constants 
take the values $k_1 = 150s^{-1}$, $k_a = 1.2 \times 10^8
s^{-1}M^{-1}$, $k_d = 10s^{-1}$ and $k_4 = 6s^{-1}$.}
\label{tab:bistab}
\end{table}

We have used an implementation of the NSM \cite{ElE}, with support for
unstructured meshes added by us. The code is written in C, and wrapped
in a MATLAB mex-file. This approach makes the definition of the
geometry, the meshing and the matrix assembly convenient.


\begin{figure}[htp]
\subfigure[t = 0s]{
\includegraphics[scale=0.55]{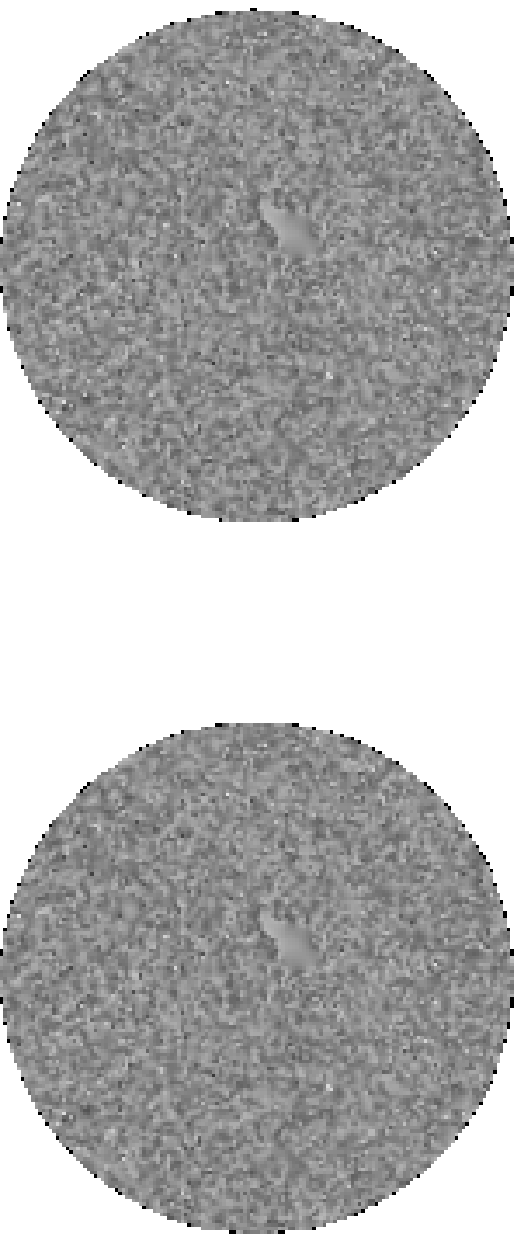}}
\subfigure[t = 2.08s]{
\includegraphics[scale=0.55]{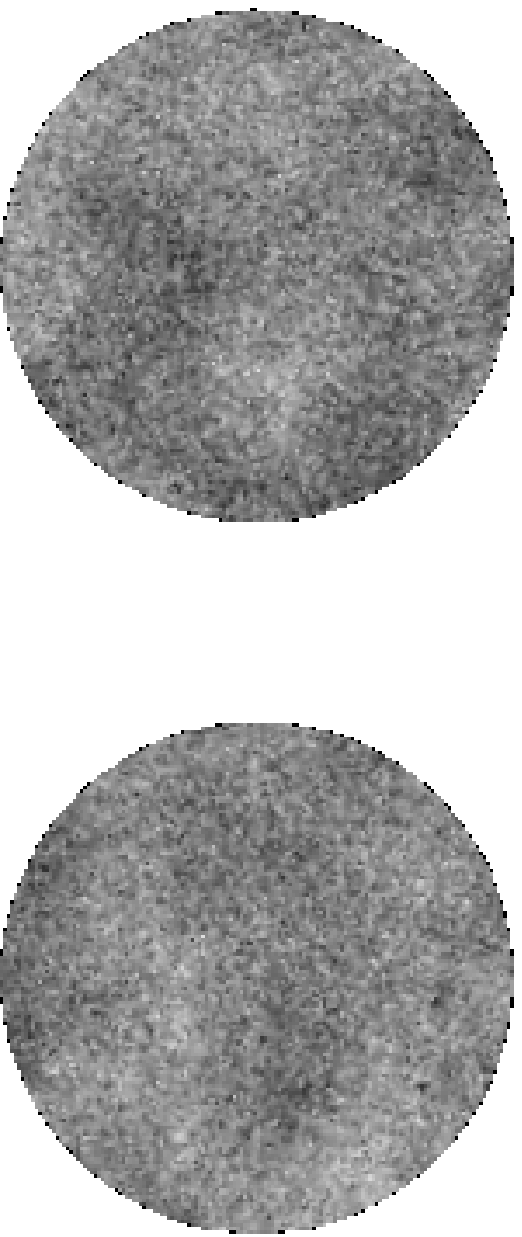}}
\subfigure[t = 8.32s]{
\includegraphics[scale=0.55]{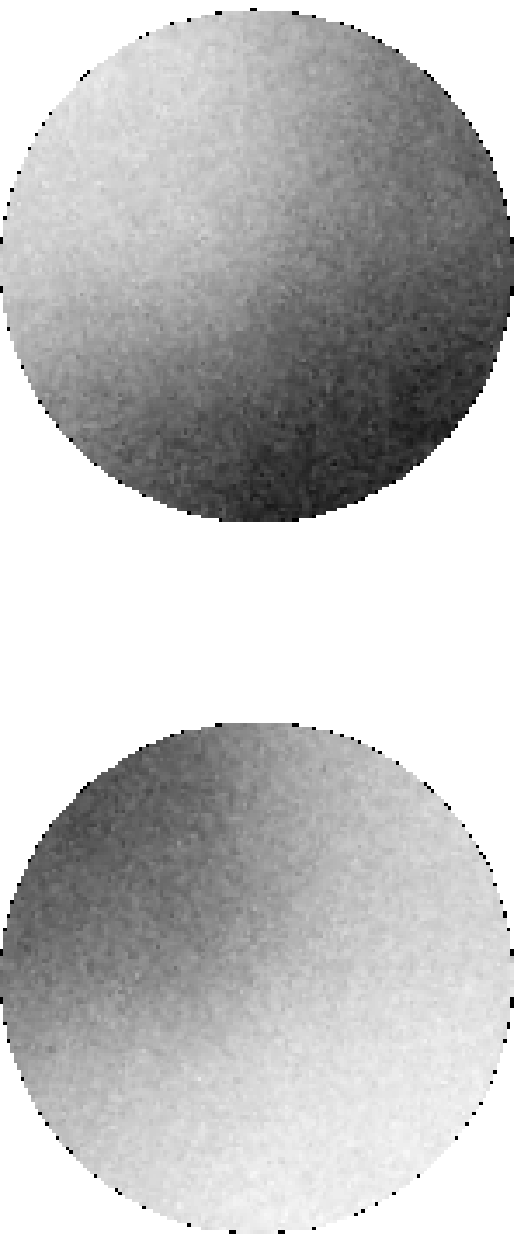}}
\caption{Snapshots of the time evolution of one trajectory with diffusion 
coefficient $\gamma = 2\times10^{-13} m^2/s$. The concentrations of the
species $A$ and $B$ are found in the upper and lower rows,
respectively. Dark areas indicate regions with a higher
concentration. Patches with the system in different phases are
formed at $t = 8.32 s$.}
\label{fig:circ1}
\end{figure}

The unstructured mesh has $K = 8849$ nodes, giving a minimal dual cell
area of $8.13 \times 10^{-16} m^2$. The mesh quality is not perfect; a
few off diagonal elements in the stiffness matrix fail to be
non-negative, giving slightly wrong diffusion rates locally. This,
however, does not seem to have any profound impact on the behavior of
the simulated system when compared to simulations on structured
triangular meshes where all coefficients are of the correct sign. The
boundaries are reflecting, corresponding to a Neumann boundary
condition in the finite element formulation of the macroscopic
equation.

The time evolution of the system is simulated on a circle with radius
$3\times 10^{-6} m$ in Figure~\ref{fig:circ1} with $\gamma =
2\times10^{-13} m^2/s$. There is a domain separation with many $A$
molecules in the lower right part and many $B$ molecules in the upper
left part.

\begin{figure}[htp] 
\subfigure[t = 0s]{
\includegraphics[scale=0.55]{fig/circ1_2.eps}}
\subfigure[t = 2.16s]{
\includegraphics[scale=0.55]{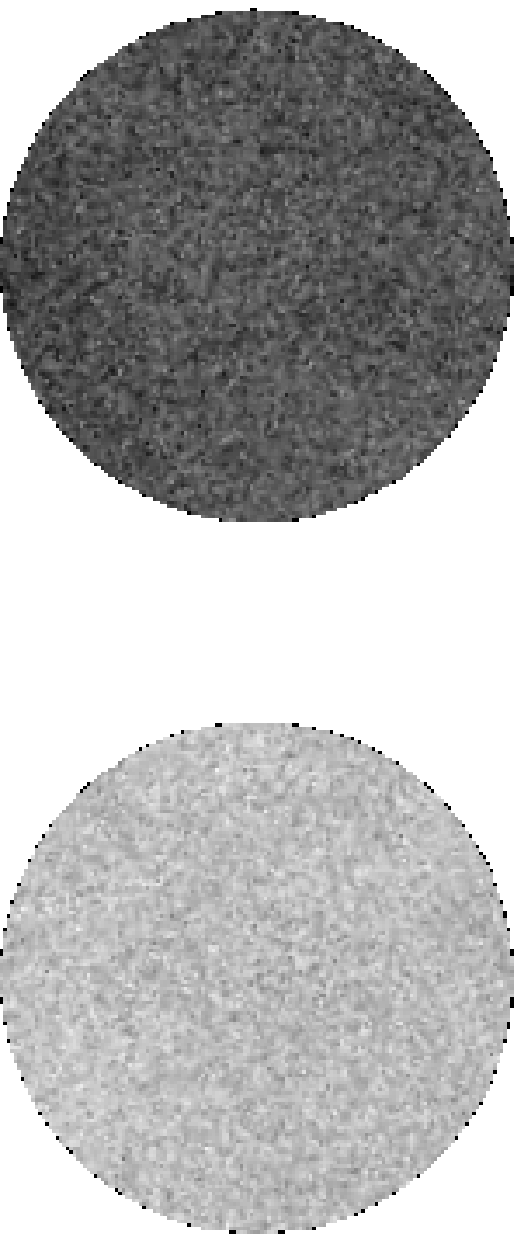}}
\subfigure[t = 7.62s]{
\includegraphics[scale=0.55]{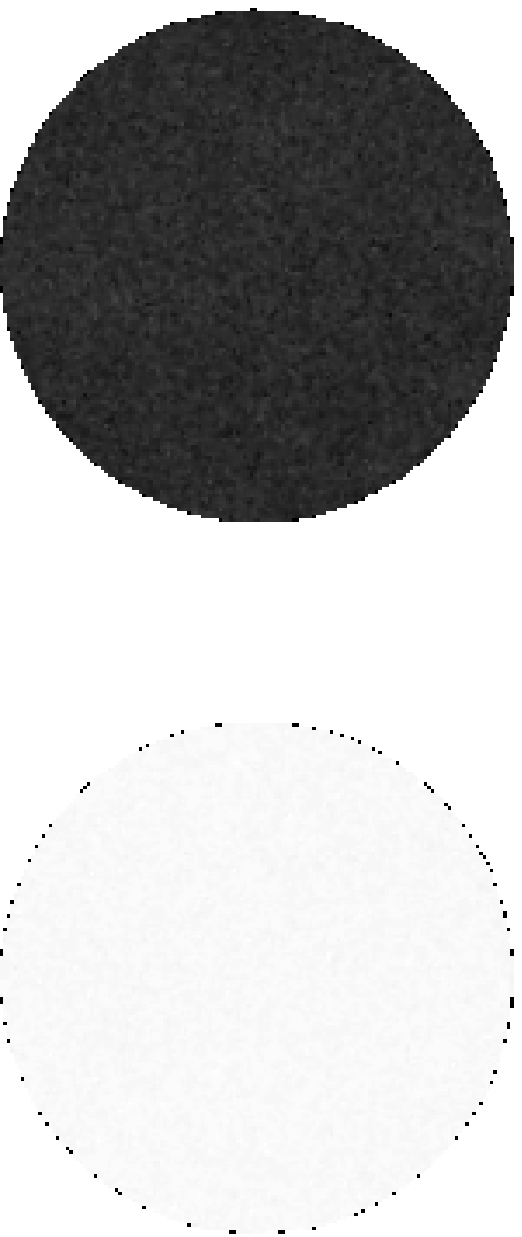}}
\caption{The same system as in Figure~\ref{fig:circ1} is simulated with 
$\gamma = 1 \times 10^{-12} m^2/s$. Here the system has not separated
into domains with different phases.}
\label{fig:circ2}
\end{figure}

In Figure~\ref{fig:circ2} the same system is simulated with fast
diffusion, $\gamma = 1\times 10^{-12} m^2/s$, and in this case the
system does not separate into domains with different phases. At the
end time of this simulation, the system is in a state where $A$
molecules dominate, and the system behaves much like in the
homogeneous case.

\subsection{The hybrid method -- metabolites and enzymes}
\label{sec:hybridres}

A model of a biochemical network with two metabolites $A$ and $B$ and
two enzymes $E_A$ and $E_B$ from \cite{fokkerplanck5} is simulated in
this example. The domain $\Omega$ is a disc with radius $\rho_{max} =
\pi^{-1/2} \approx 0.564$ and area $|\Omega| = 1$. The reactions
are summarized in Table~\ref{tab:chemreac}. Initially, the
concentrations $a$ and $b$ are constant and the enzyme concentrations,
$e_{A}$ and $e_{B}$, are zero in every cell.

\begin{table}[H]
\begin{align*}
&\begin{array}{lllll}
  A \overset{\mu}{\underset{w_{1}}{\rightleftharpoons}} \emptyset &
  B \overset{\mu}{\underset{w_{2}}{\rightleftharpoons}} \emptyset & 
  A+B \xrightarrow{k_{2}} \emptyset &
  E_A \overset{\mu}{\underset{w_{3}}{\rightleftharpoons}} \emptyset &
  E_B \overset{\mu}{\underset{w_{4}}{\rightleftharpoons}} \emptyset
\end{array}                                                             \\
  w_{1}(a,e_A) &= k_A e_A/(1+a k_I^{-1}) \qquad
  w_{3}(a)  = H(0.2-\rho) k_{E_A}/(1+a k_R^{-1})                        \\
  w_{2}(b,e_B) &= k_B e_B/(1+b k_I^{-1}) \qquad
  w_{4}(b)  = H(\rho-0.4) k_{E_B}/(1+b k_R^{-1})
\end{align*}
\caption{Reaction channels for the network. The concentrations of the 
species $A, B, E_A,$ and $E_B$ in a cell are $a, b, e_A, e_B$. The
reaction constants are $k_A = k_B = 3\zeta$, $\mu = 0.002\zeta$, $k_2
= 0.0005\zeta^{2}$, $k_{E_A} = k_{E_B} = 0.5\zeta$, $k_I = 60/\zeta$,
$k_R = 30/\zeta$ and the diffusion constant is $\gamma = 10^{-4}$. The
domain is covered by $K = 80$ dual cells and the scaling with the
average size of a cell $\zeta \equiv |\Omega|/K$ is done to define the
unit scale of the problem. $H$ denotes the Heaviside function; the
enzyme $E_A$ is thus produced only in the center of $\Omega$, $0 \le
\rho < 0.2$, and $E_B$ is created only close to the boundary, $0.4 \le
\rho \le \rho_{max}$.}
\label{tab:chemreac}
\end{table}
 
\begin{figure}[htp]
  \centering
  \subfigure{\includegraphics[width = 5cm]{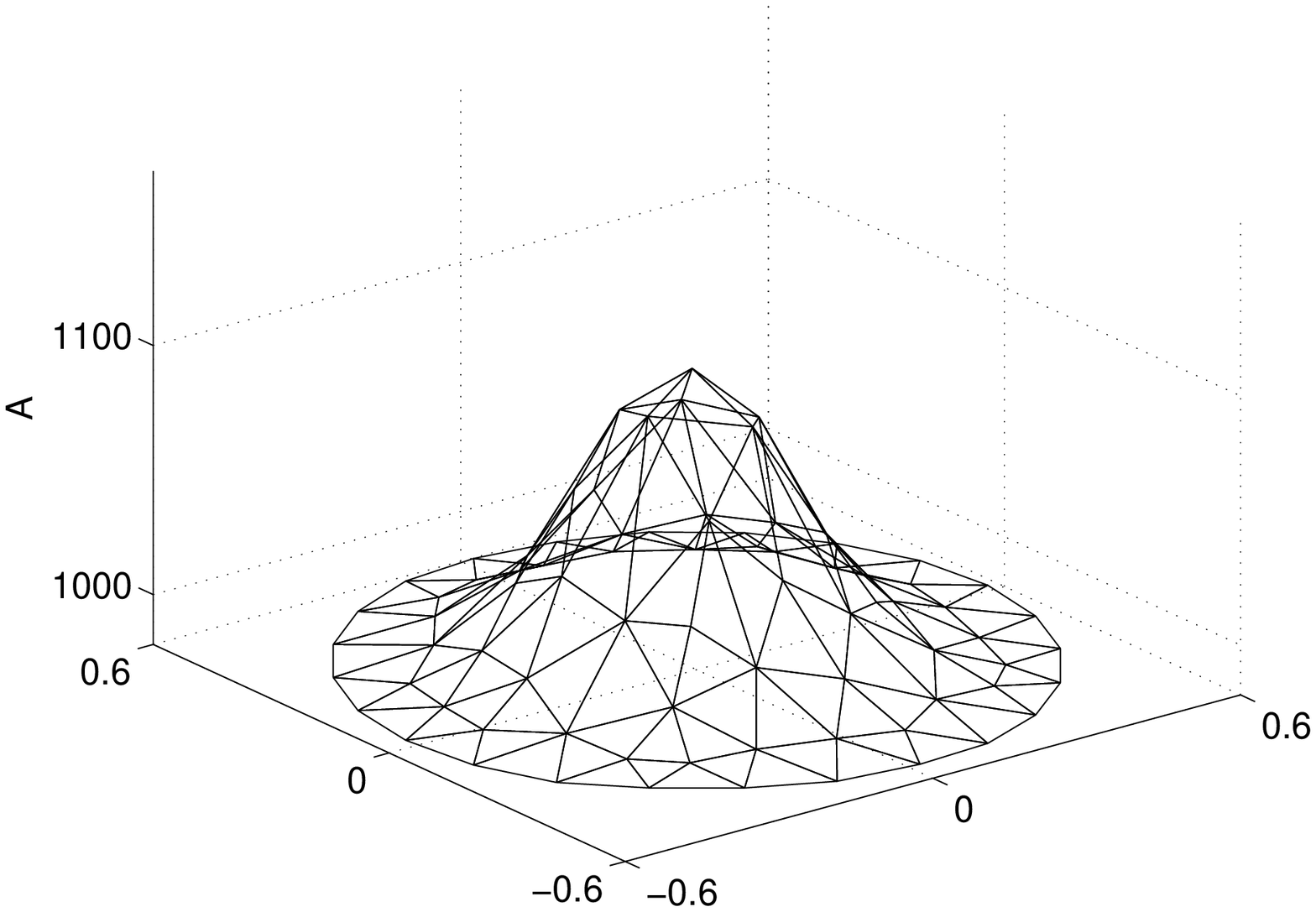}}
  \subfigure{\includegraphics[width = 5cm]{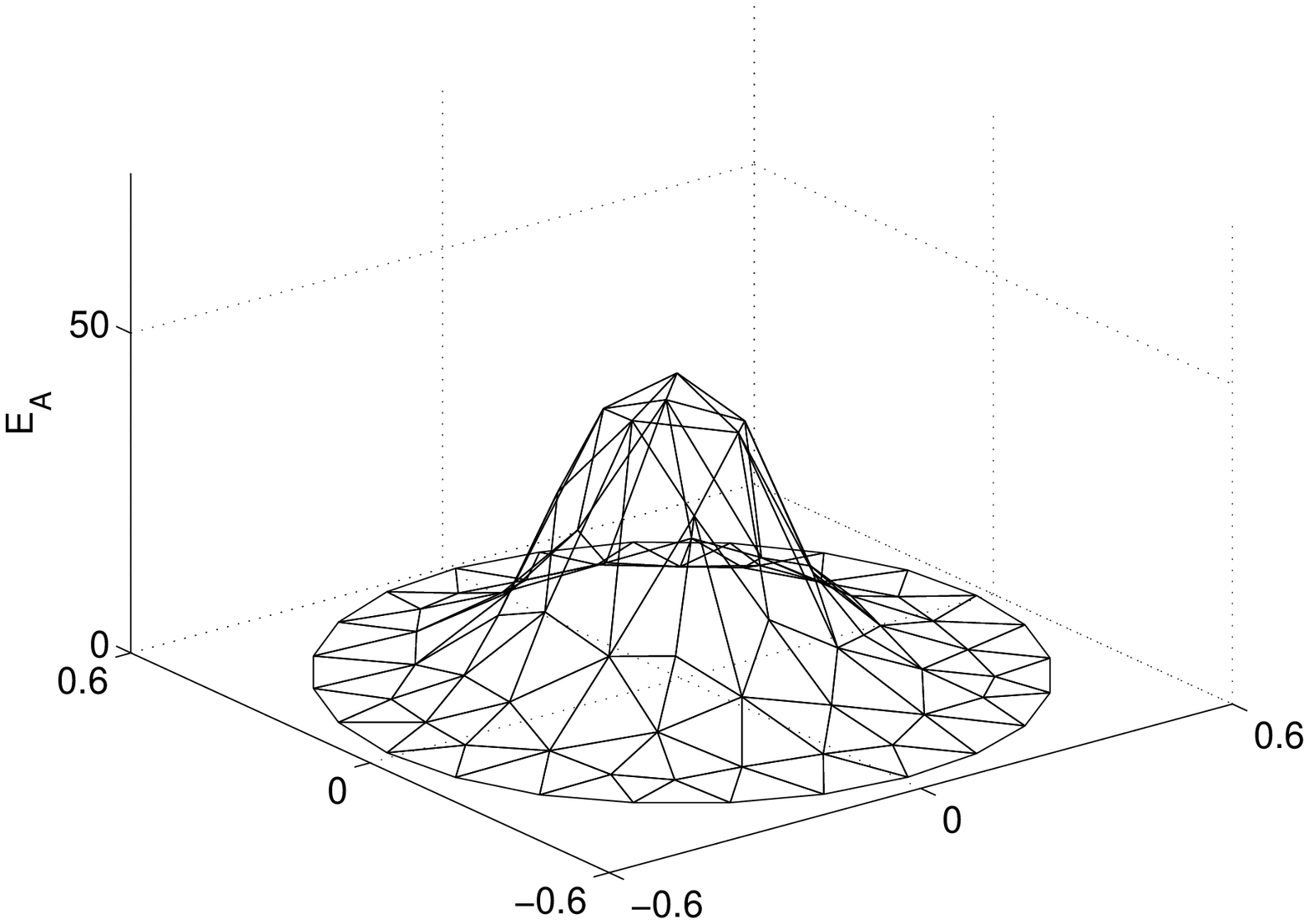}}
  \\
  \subfigure{\includegraphics[width = 5cm]{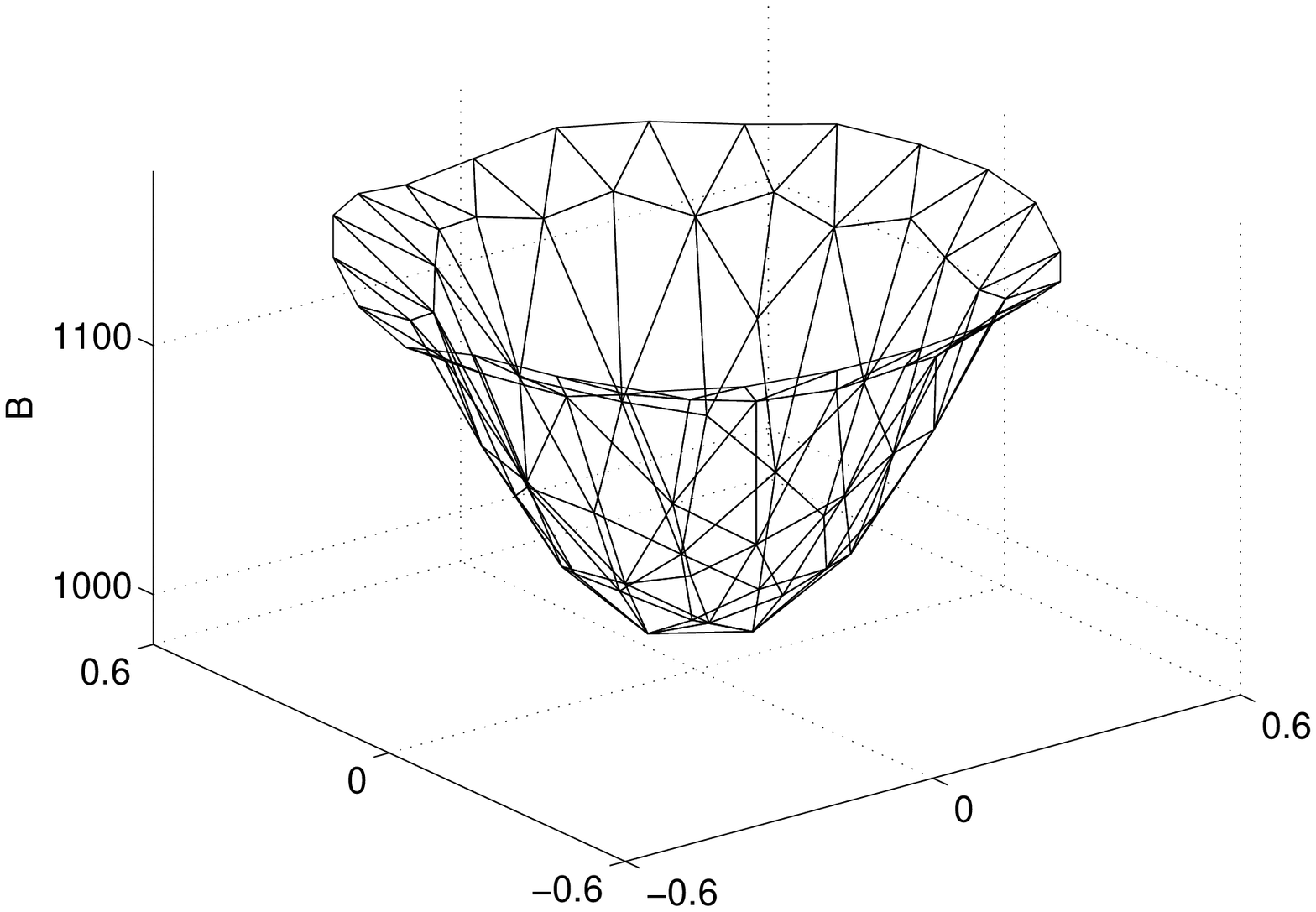}}
  \subfigure{\includegraphics[width = 5cm]{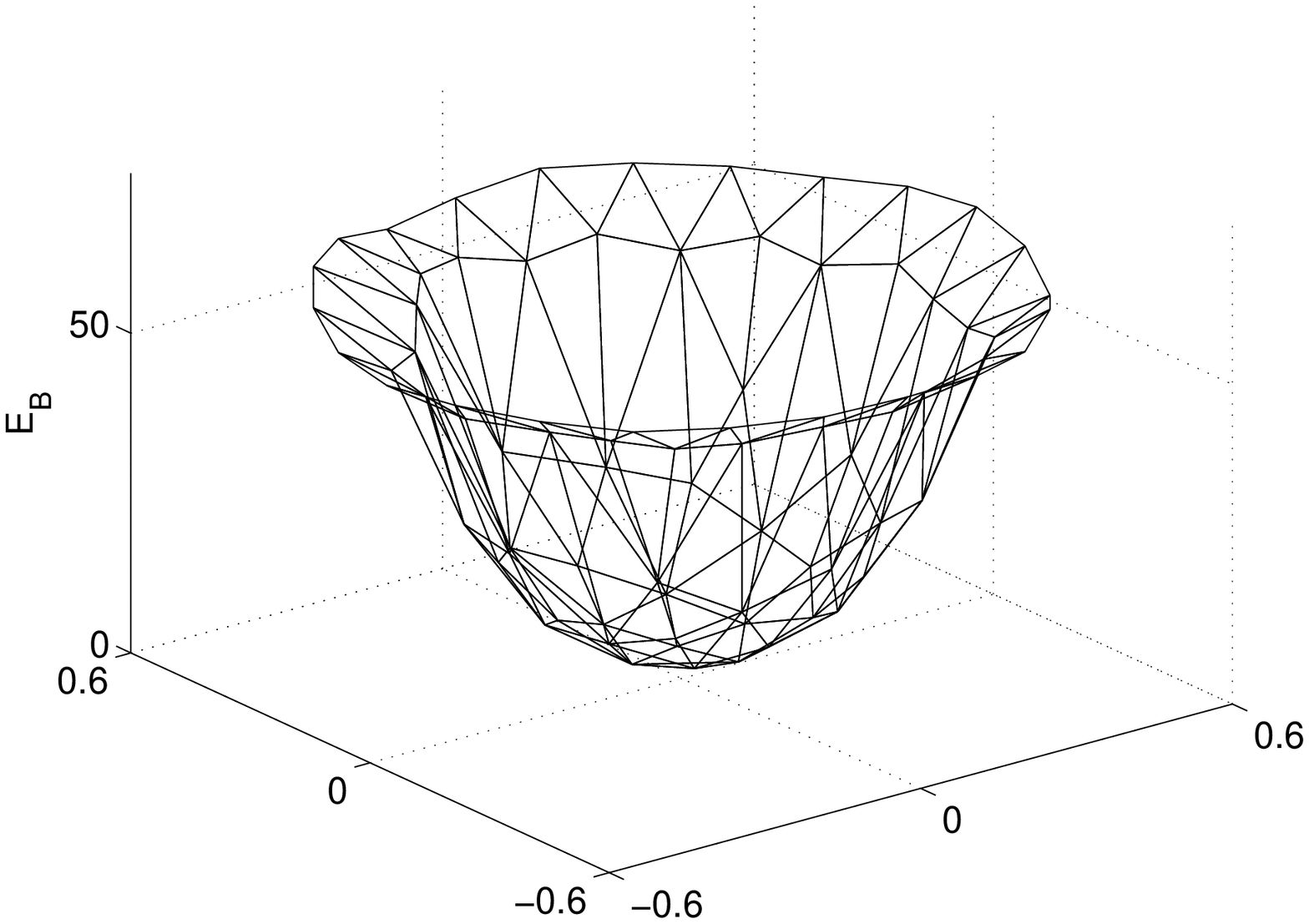}}
  \caption{The concentrations of the metabolites $A$ (upper left) and
  $B$ (lower left) and enzymes $E_A$ (upper right) and $E_B$ (lower
  right) at $t = 200$.}
  \label{fig:AB}
\end{figure}

The chemical reactions and the diffusion of the enzymes are simulated
with SSA and the diffusion of $A$ and $B$ is modeled by the diffusion
PDE in a straightforward MATLAB implementation of the hybrid
method. The effect of the enzymes in different parts of $\Omega$ is
demonstrated in Figure~\ref{fig:AB} with $M = 10^4$ realizations.

A stochastic reference solution $\fatphi_{i \cdot}^s,\; i = 1,2,3,4,$
averaged over $M = 10^4$ trajectories is integrated directly to $t =
200$ with SSA. Hybrid solutions $\fatphi_{i \cdot}^{\Delta t}(200)$
are computed for different $\Delta t$ with $M = 10^4$ for comparison
with $\fatphi_{i \cdot}^s$. Table~\ref{tab:diff2} displays the
differences between the stochastic and the hybrid solutions.

\begin{table}[htp]
  \centering
  \begin{tabular}{|c|c|c|c|c|c|c|}\hline
    $\Delta t$&0.1&1&5&20&40&100\\ \hline
    $\delta_t$&0.024&0.024&0.024&0.024&0.025&0.030\\ \hline
  \end{tabular}
  \caption{The difference between the stochastic solution and the
  hybrid solution in the $\ell_2$-norm, $\delta_t = \max_i
  {\|\fatphi_{i \cdot}^s-\fatphi_{i \cdot}^{\Delta t} \|_2}/ (\max_{j}
  \fatphi_{ij}^s-\min_{j}\fatphi_{ij}^s)$ at $t = 200$ for different
  $\Delta t$. The stochastic errors dominate in $\delta_t$ for $\Delta
  t\lesssim 40$.}
  \label{tab:diff2}
\end{table}

\begin{figure}[htp]
  \centering
  \subfigure[]{\includegraphics[width = 5.5cm]{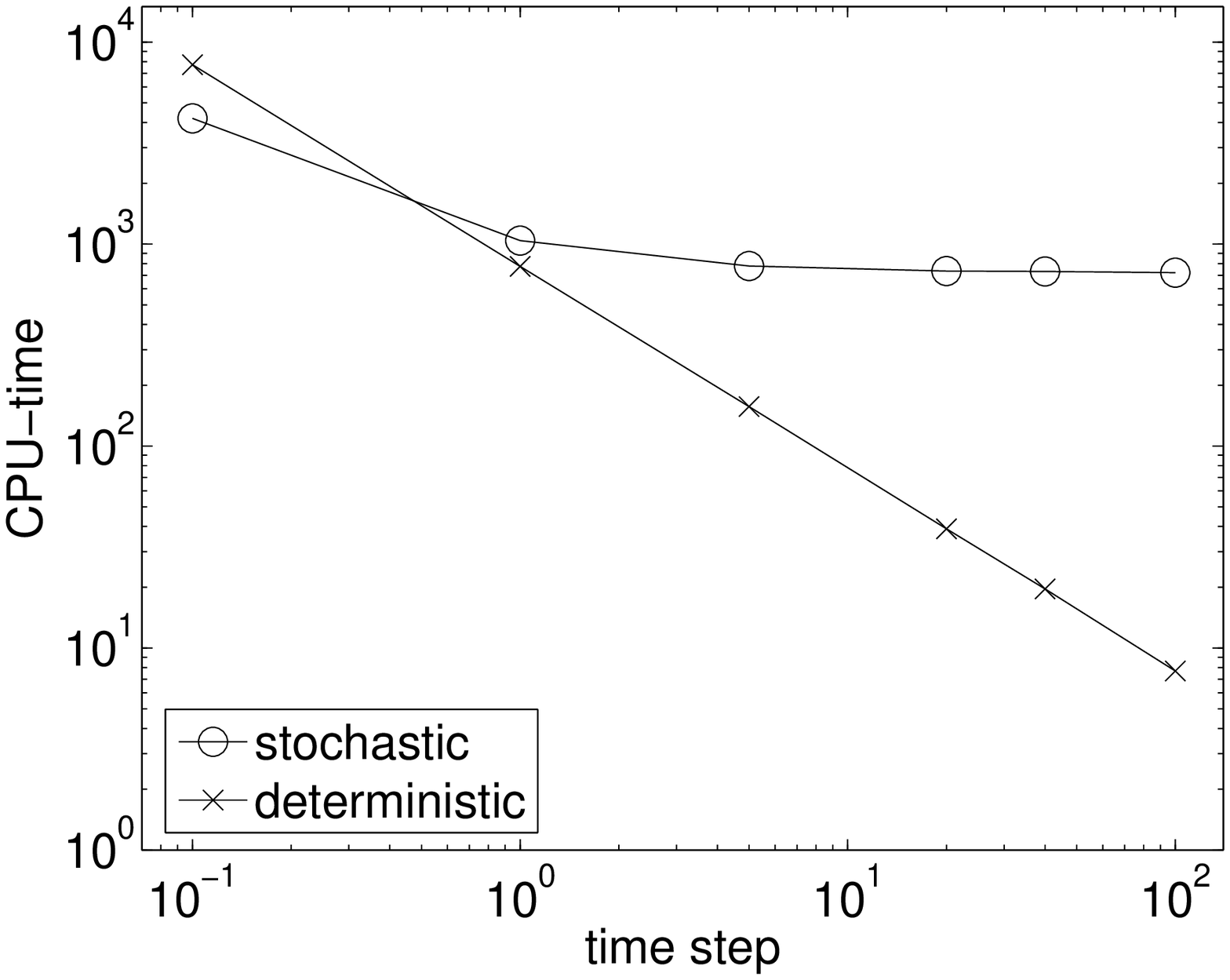}}
  \subfigure[]{\includegraphics[width = 5.5cm]{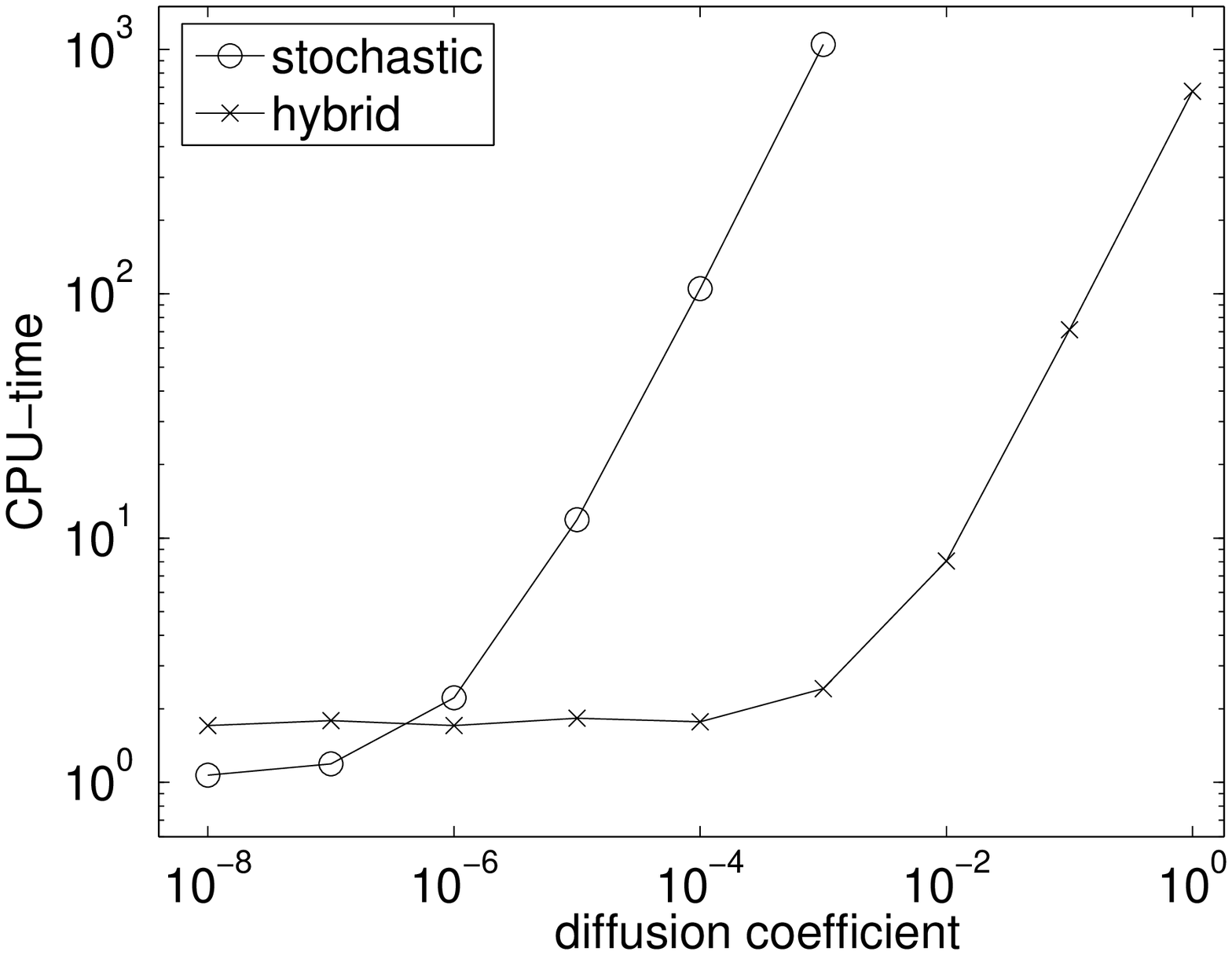}}
  \caption{Run times versus the time step for the stochastic and
  deterministic part of the hybrid algorithm (a) and run times versus
  the diffusion coefficient $\gamma$ for SSA and the hybrid algorithm (b).}
  \label{fig:runtimes}
\end{figure}

The computational work for the stochastic and the deterministic parts
of the algorithm is compared in Figure~\ref{fig:runtimes}~(a) at
$t=200$ with $\gamma = 10^{-4}$ and $M = 10^4$.  The work in the
deterministic part is less than the work in the stochastic part when
$\Delta t>1$, and decreases as $\Delta t^{-1}$ when $\Delta t$
increases. The work for the stochastic part tends to a limit since the
extra effort for restarting SSA in each step becomes negligible when
the time step is longer.

Figure~\ref{fig:runtimes}~(b) displays the total CPU-time for
stochastic and hybrid solutions for different diffusion coefficients
$\gamma$ at $t = 10$ with $M = 10^3$. We use $\Delta t = 5$ for the
hybrid algorithm, i.e. a time step well below the upper limit for
which the stochastic errors dominate for $\gamma = 10^{-4}$.  The
stochastic solutions are integrated directly to $t = 10$.  When the
diffusion of the molecules is the major part of the computational work
for larger $\gamma \gtrsim 10^{-5}$, then replacing the diffusion for
the species with large copy numbers by the macroscopic diffusion
reduces the CPU-time by up to 1000 while retaining small differences
between the solutions, see Table~\ref{tab:diff2}. The number of
molecules of the enzymes $n_e = (\fatx_{3 \cdot}+\fatx_{4
\cdot})\fate_1$ is about 0.001 of the number of metabolite molecules
$n_m=(\fatx_{1 \cdot}+\fatx_{2 \cdot})\fate_1$ at $t = 10$. This
difference explains the remarkable improvement in efficiency in
diffusion dominated regimes. From the figure and the discussion, the
CPU-times $T_{SSA}$ and $T_{hyb}$ for SSA and the hybrid method,
respectively, are approximately
\[
  T_{SSA}\approx c_{s0}+c_s\gamma (n_m+n_e),\;
  T_{hyb}\approx c_{h0}+c_s\gamma n_e,
\]   
where $c_{s0}$ and $c_{h0}$ are small and are mainly due to the
chemical reactions.
The concentration of the enzymes increases when the integration is
continued to $t = 200$. The speedup by the hybrid method with $\Delta
t = 5$ is then about 35. Since $n_m/n_e\approx 30$ this is in
agreement with the estimates of $T_{SSA}$ and $T_{hyb}$ above.

%% file: concl.tex
\section{Conclusions}
\label{sec:concl}

The RDME is discretized on an unstructured mesh for better geometric
flexibility than a Cartesian mesh. The diffusion coefficients in the
RDME are derived from a FEM discretization of the
Laplacian. Stochastic simulation can then be used to determine a
number of trajectories of the mesoscopic system. If the copy number is
large for some chemical species, then a hybrid method integrating the
diffusion part deterministically reduces the computing time
substantially, especially when the diffusion constant is large.

The method is applied to three different systems. The convergence of a
system without reactions to the macroscopic solution is shown in the
first example. In the second example, we consider a biochemical system
with both diffusion and reactions. We illustrate that our method can
be efficiently applied to a previously studied bi-stable system, and
that the results obtained with our code are in line with those
computed on structured meshes with the freely available software
MesoRD \cite{HFE}. In a final example, the performance of the hybrid
method is compared to SSA for a system with four species. The hybrid
algorithm is up to three orders of magnitude faster.

In this paper we have considered examples in two space dimensions
only. Realistic modeling of the reaction networks in e.g.~bacteria
typically require 3D simulations. With our approach, the extension to
3D is straightforward and will be reported in a forthcoming paper. The
method relies on the FEM discretization of the diffusion equation and
a variety of existing software can be used to specify the geometry,
construct the mesh and postprocess the result. Presently, we handle
diffusion with a uniform diffusion constant but there is no
complication in considering space-dependent diffusion or adding
convection or to let the diffusion be different for different species.

In the spatially homogeneous case, stiffness arises from the presence
of fast chemical reactions. For the RDME, the stiffness in addition
increases with the resolution of the mesh. We have proposed a hybrid
method to reduce simulation time when some species are present in
large copy numbers. For the homogeneous case, many approximative
schemes have been developed for systems with time scale
separation. For spatially dependent systems, no sharp separation in
slow and fast events can generally be made. Instead we rather have a
continuum of scales and multiscale simulation techniques have to be
developed anew.